%% file: mk-10_paper.tex
\documentclass[11pt]{article}
\long\def\kmcomment#1{}
\input{./kmme11.tex}

\input{./miz-pers.mac}
\usepackage{graphicx,color,verbatim,setspace,supertabular,ntheorem} 
\usepackage{latexsym,amssymb,amsmath,alltt}
\usepackage{url}
\newcommand{\myd}{d\, }
\renewcommand{\tz}[2]{\widehat{z}^{#1}_{#2}}
\newcommand{\nz}[2]{{z}^{#1}_{#2}}
\newcommand{\kmqed}{\hfill \rule{1ex}{1.5ex}\par}

{
\theorembodyfont{\rmfamily}%{\bfseries}{\rmfamily}{\sffamily}
\theoremstyle{plain}%{break}  
\newtheorem{thm}{Theorem}[section]

\newtheorem{exam}[thm]{Example}
\newtheorem{defn}[thm]{Definition}
\newtheorem{kmProp}[thm]{Proposition}

\newtheorem{myRemark}[thm]{\textbf{Remark}}
}

\newcommand{\CGFF}[2]{\text{C}^{#1}_{\rm GF}({\kham}_2^0,{\fraksp}(2,\mR))_{#2}}
\newcommand{\HGFF}[2]{\text{H}^{#1}_{\rm GF}({\kham}_2^0,{\fraksp}(2,\mR))_{#2}}

\newcommand{\CGF}[4]{\text{C}^{#1}_{\rm GF}({\kham}_{#2}^{#3},{\fraksp}(#2,\mR))_{#4}}
\newcommand{\HGF}[4]{\text{H}^{#1}_{\rm GF}({\kham}_{#2}^{#3},{\fraksp}(#2,\mR))_{#4}}

\newcommand{\cgf}[4]{\text{C}^{#1}_{\rm GF}({\kham}_{#2}^{#3})_{#4}}
\newcommand{\hgf}[4]{\text{H}^{#1}_{\rm GF}({\kham}_{#2}^{#2})_{#4}}

\title{Higher weight Gel'fand-Kalinin-Fuks classes of formal Hamiltonian 
vector fields of symplectic $\mR^2$}

\author{
Kentaro Mikami\thanks{ 
  Department of Computer Science and Engineering 
  Akita University, partially supported by Grant-in-Aid for 
  Scientific Research (C) of JSPS, No.23540067 and No.20540059}
\and   
Yasuharu Nakae\thanks{ 
  Department of Computer Science and Engineering 
  Akita University} 
\and
Hiroki Kodama\thanks{
Graduate School of Mathematical Sciences, the University of Tokyo,
partially supported by ``FIRST Program'' of JSPS.
}
  }

\date{}

\begin{document}
\maketitle
\renewcommand{\thefootnote}{\fnsymbol{footnote}}
\begin{abstract}
In ``The {G}el'fand-{K}alinin-{F}uks class and characteristic classes of
transversely symplectic foliations'', arXiv:0910.3414, (October 2009) by
D.~Kotschick and S.~Morita, the relative {G}el'fand-{K}alinin-{F}uks cohomology
groups of the formal Hamiltonian vector fields without constant vector fields
on 2n-plane were characterized by two parameters, one is degree and the other
is weight. And they obtained those cohomology groups of the 2-plane while their weight $\le$ 10.  

In this paper, for those cohomology groups of the 2-plane, we succeeded in
determining the dimension of cochain complexes by Sp(2,R)-representation theory
for their weight even less than 50, thus, we manipulate the Euler
characteristic numbers.   We also decide our relative
{G}el'fand-{K}alinin-{F}uks cohomology groups until whose weight $<$ 20 by
getting a concrete matrix representation of the coboundary operator.
\end{abstract}

\section{Introduction}% 
Since 2000, many works have been achieved
on Gel'fand-Kalinin-Fuks cohomology.
S.~Metoki in his doctoral thesis \cite{metoki:shinya}
found a new exotic class, which is now called the Metoki class. 
M.~Takamura \cite{M:Takamura} showed that 
the relative cohomology of 
Lie algebra of formal contact vector fields with
respect to formal Poisson vector fields is trivial.    
In the paper of D.~Kotschick and S.~Morita \cite{KOT:MORITA},
the relative Gel'fand-Fuks cohomology groups
were characterized by two parameters; the degree and the weight, 
and they obtained those cohomology groups of the $2$-plane with
weight $\leq 10$.
We are ambitious to compute more higher weight cases.

In this paper, we recall fundamental facts on Hamiltonian formalism;
particularly we review group actions very carefully.
And then we concentrate
on the weight. We show that the weight corresponds to Young diagrams,
and we construct the generating function for the weight. 
We make use of representation theory of $Sp(2n)$, 
especially the fact that
the irreducible representations are parameterized 
by the Young diagrams of height at most $n$. 
In the final section we exhibit our main result,
which we obtain by means of computer calculation.

\medskip

\texttt{Points of this revised version (on Feb 2014):}

(1) Fixing notations of Lie algebras of formal Hamiltonian
vector fields: 
\begin{center}
\begin{tabular}{|c|ccc|} 
        \hline
this article & $\ds\kham_{2n}     $ & $\ds\kham_{2n}^{0}$& $\ds\kham_{2n}^{1}$
        \\ \hline older & $\ds\kham_{2n}^{0}$ & $\ds\kham_{2n}^{1}$& \\ \hline 
                        & $\oplus_{\ell=1}^{\infty} S^{\ell}$   
                        & $\oplus_{\ell=2}^{\infty} S^{\ell}$   
                        & $\oplus_{\ell=3}^{\infty} S^{\ell}$ \\\hline
\end{tabular}
\end{center} 
where $\ds S^{\ell}$ means the set of $\ell$-th homogeneous polynomials, to
\cite{Kont:RW} or \cite{KOT:MORITA}.  
\kmcomment{
M. Kontsevich, {\em Rozansky-Witten invariants via formal geometry},
Composito Math. 115 (1999), 115--127.  
}

(2) An announcement of Betti numbers for the weight 20. 
\section{Preliminaries}
\subsection{Lie algebra
Cohomology and 
Gel'fand-Kalinin-Fuks 
cohomology}  
Take a Lie algebra $\mathfrak g$ over $\mR$.
Let $(\rho,W)$ be a representation of $\mathk g$. Namely, $\rho$ is a Lie
algebra homomorphism of $\mathk g$ into the Lie algebra $\text{End}(W)$.   
For each $k\in \mZ_{\ge 0}$, $$\ds  \text{C}^{k}( {\frakg}) := \{ \sigma :
\underbrace{
{\mathk g} \times \cdots \times {\mathfrak g}}_{k-\text{times}} \rightarrow W \mid
\text{alternative
and $\mR$-multilinear}\}. $$  For each k-th cochain 
$\sigma \in \text{C}^{k}( {\frakg}) $, we define 
\begin{align*}
 (\myd \sigma)(X_0, \ldots, X_k) :=& \sum_{i=0}^k (-1)^i \rho(X_i) \sigma (\ldots
\widehat{X_i}\ldots)\\
& + \sum_{i<j} (-1)^{i+j} \sigma([X_i,X_j]\ldots
\widehat{X_i} \ldots \widehat{X_j} \ldots) 
\end{align*}
it is known that $\myd$ satisfies $\ds  \myd^2=0 $ and defines the
cohomology groups of a Lie algebra 
$\mathk g$ with respect to $(\rho, W)$.  
In this paper, hereafter we only deal with  
the trivial representation of Lie algebra, i.e., $W=\mR$ and
$\rho=0$.  

Let $\mathk k$ be a subalgebra of $\mathfrak g$.  Define 
$$\ds  \text{C}^{m}( {\mathk g}, {\mathfrak k} ) := \{ \sigma \in
C^{m}({\mathk g}) \mid i_X\sigma =0, i_X \myd \sigma =0\quad
(\forall  X\in {\mathk k})\}$$ and we get the relative cohomology groups $\ds 
\text{H}^{m}({\frakg}, {\frakk})$.  
Let $K$ be a Lie group of $\mathk k$. Then we also consider 
$$\ds  \text{C}^{m}( {\mathk g},  K ) := \{ \sigma \in
        \text{C}^{m}({\mathk g}) \mid i_X\sigma =0\  
(\forall  X\in {\mathk k}), Ad_k^{*} \sigma = \sigma \ (\forall k\in K)
\}$$ and we get the relative cohomology groups $\ds 
\text{H}^{m}({\mathk g}, K)$. If $K$ is connected, those are identical.  
If $K$ is a closed subgroup of $G$, then $\text{C}^{\bullet}({\mathk g}, K)
= \Lambda ^{\bullet} ( G/K) ^{G} $ (the exterior algebra of
$G$-invariant differential forms on $G/K$).

Take a differentiable manifold $M$ and consider the space ${\mathk X}(M)$ 
of vector fields of $M$. Then  ${\mathk X}(M)$ forms a Lie algebra by
Jacobi-Lie bracket.  Thus, we can consider the Lie algebra cohomology of
${\mathk X}(M)$.  But, the cochain complex is huge, so we add some
restriction, ``continuity'' by $C^{\infty}$-topology.  
Let ${\mathk g}$ be a subalgebra of ${\mathfrak X}(M)$.   
For each $k\in \mZ_{\ge 0}$, $$\ds  \text{C}^{k}( {\frakg}) := \{ \sigma :
\underbrace{
{\mathk g} \times \cdots \times {\mathfrak g}}_{k-\text{times}}
\rightarrow \mR  \mid
\text{continuous, alternative
and $\mR$-multilinear}\}. $$  For each k-th cochain 
$\sigma \in \text{C}^{k}( {\frakg}) $, we define 
\begin{align*}
 (\myd \sigma)(X_0, \ldots, X_k) :=&   \sum_{i<j} (-1)^{i+j} \sigma([X_i,X_j]\ldots
\widehat{X_i} \ldots \widehat{X_j} \ldots) .
\end{align*}
It is known that $\myd$ satisfies $\ds  \myd^2=0 $ and defines the
cohomology group, called Gel'fand-Kalinin-Fuks cohomology of ${\mathk
X}(M)$.  There are also relative versions.

\subsection{Recall of Hamilton formalism}

Let $(M,\omega)$ be a symplectic manifold, namely, $\omega$ is a
non-degenerate closed 2-form on $M$, and so $\dim M$ is even.  
We denote the group of symplectic automorphisms of 
$(M,\omega)$ by $Aut(M,\omega)$.  
By
$\ds  {\frakaut}(M,\omega)$, we denote the space of vector
fields on $M$ satisfying $\ds  {\cL}_X \omega=0$ (infinitesimal
automorphism of $\omega$).  

Since $\ds  {\cL}_{[X,Y]} = [{\cL}_X, {\cL}_Y]$ holds, $\ds 
{\frakaut}(M,\omega)$ forms a Lie algebra.  

For $\forall f$ on $M$, we
have the Hamiltonian vector field on $(M,\omega)$ defined by
$\ds  \omega(\Hd _f,\cdot) = df$, and the Poisson bracket given by
$\ds  \Pkt{f}{g}:= \omega( \Hd _f, \Hd _g)$.  Since $\ds 
\Hd _{\Pkt{f}{g}} = - [ \Hd _f, \Hd _g]$ holds,  
the Hamiltonian vector fields of
$(M,\omega)$  form a Lie subalgebra of $\ds 
{\frakaut}(M,\omega)$. The correspondence $\ds  f \mapsto
-\Hd_f$ is a Lie algebra homomorphism and the kernel is $\mR$ when $M$ is
connected.  

It holds   
$ \ds  \varphi  \Hd _f = \Hd _ { f \circ \varphi^{-1}}$ for
each $\varphi  \in Aut(M,\omega)$, and   
$ \ds  \Pkt{f}{g} \circ \varphi  = 
\Pkt{f\circ \varphi}
{g\circ \varphi}$ for $f,g\in \text{C}^{\infty}(M)$.

Let $K$ be a Lie subalgebra of $Aut(M,\omega)$ with 
its Lie algebra $\ds  \frakk$.  
For each  $\xi\in {\frakk}$, the
fundamental vector field on $M$, say $\ds  {\xi}_{M}$, is defined
by $\ds   \xi_{M} := \frac{d}{dt} \exp ( t \xi)_{M \mid t=0}$.  
They satisfy $\ds  {\cL}_{ \xi_{M}} \omega = 0$ and $\ds 
[\xi,\eta]_{M} = - [\xi_{M}, \eta_{M}] $, thus  form  a subalgebra of
$\ds  {\frakaut}(M,\omega)$.  The momentum mapping $J$ (if
exists) is a
map from $\ds  M \rightarrow \frakk^{*}$
satisfying $$ d \hat{J} (\xi ) = \omega( \xi_{M}, \cdot),
\quad\text{i.e.,}\quad \xi_M = \Hd _{\hat{J}(\xi)}$$
where 
$
\ds  \hat{J}(\xi) $ is defined as  
$\ds  \langle \hat{J}(\xi), m\rangle := \langle \xi,
J(m)\rangle$ for $\xi\in \frakk, m \in M$.  
The above definition 
 means $J$ provides with a Hamiltonian potential for each fundamental vector field
of $K$ and 
$
\ds  \hat{J}(\xi) $ is defined as  
$\ds  \langle \hat{J}(\xi), m\rangle := \langle \xi,
J(m)\rangle$ for $\xi\in \frakk, m \in M$.  
If 
$\ds  J$ is $K$-equivariant, i.e., 
$ J( a\cdot m) = Ad_{a^{-1}}^{*} ( J(m) )$ for $\forall a\in K, m\in M$, 
then $\ds  
\hat{J}([\xi, \eta]) = \Pkt{ \hat{J}(\xi)}{\hat{J}(\eta)}$ 
hold for $\xi,\eta\in \frakk$, and vice versa if $K$ is connected.

In local, by Darboux's theorem we always have a local coordinates
$\ds  q^1,\ldots, q^n, p_1,\ldots, p_n$ such that 
$\ds  \omega( \frac{\pdel}{ \pdel q^i}, \frac{\pdel }{\pdel
p_i}) = 
- \omega( \frac{\pdel }{\pdel p_i}, 
\frac{\pdel }{\pdel q^i}) = 1$, the others are 0 
 and so the Hamiltonian vector
field is given by  
$$\ds  \Hd _f = \sum_{i=1}^n \left( \frac{\pdel
f}{\pdel p_i} \frac{\pdel}{\pdel q^i} - \frac{\pdel f}{\pdel q^i}
\frac{\pdel}{\pdel p_i}\right)$$ 
and the Poisson bracket is given by
$$\ds  \Pkt{f}{g} 
= \sum_{i=1}^n \dfrac{ \pdel ( f,g )}{ \pdel ( q^i, p_i )} 
= \sum_{i=1}^n \left( \frac{\pdel f}{\pdel
q^i} \frac{\pdel g}{\pdel p_i} - \frac{\pdel f}{\pdel p_i} \frac{\pdel
g}{\pdel q^i}  \right) .
$$

If $\ds  M= 
(\mR^{2n},\text{linear symplectic structure})$, the space of the
Hamiltonian vector fields of $(M,\omega)$  coincides with $\ds 
{\frakaut}(M,\omega)$, because of the space is connected and
1-connected.   
Let $K$ be the linear symplectic group of $\ds  (\mR^{2n},
\omega)$ with its Lie algebra $\ds  \frakk$, i.e., 
$\ds  K = Sp(2n,\mR)$ and 
$\ds  \frakk  = \fraksp(2n,\mR)$.   
The equivariant (co-)momentum mapping is given by  
$$
\hat{J}(\xi) (q,p) = \frac{-1}{2} (q,p) \left(
\text{matrix representation of } \omega \right) 
\xi
\begin{pmatrix} q \\ p\end{pmatrix}
$$ here $q$ is the natural coordinate of $\ds  \mR^n$, 
$\ds  \hat{J}$ is a Lie algebra
homomorphism from the Lie algebra $\fraksp(2n,\mR)$ into $\ds 
C^{\infty}(M),$ with the Poisson bracket.     
The Hamilton potential of Hamiltonian vector field $\ds 
[\xi_{M}, \Hd _f]$ is given by $\ds  - \Pkt{\hat{J}(\xi)}{f}$,
because of $\ds  [\xi_{M}, H_f] = [\Hd_{\hat{J}(\xi)}, \Hd_f] 
= - \Hd_{ \Pkt{ \hat{J}(\xi)}{f} } 
$.

\section{Cochain complexes, weight and relative cochains} 
\subsection{Cochain complexes}
When $M=\mR^{2n}$,  
the Hamiltonian potential for each Hamiltonian vector field 
is unique up to
constant,  we consider the Lie subalgebra 
$\ds {\frakham_{2n}}$ of 
Hamiltonian vector fields 
is isomorphic with the formal polynomial space. Then the space is 
Lie algebra isomorphic with the space of 
$\ds  
\mR[[q^1,\ldots, q^n, p_1,\ldots, p_n]]/\mR$ quotiented by $\mR$,   
where the Lie bracket is given by the Poisson bracket.

We see that 
$$\ds   \mR[[q^1,\ldots, q^n, p_1,\ldots, p_n]]/\mR =
\text{the completion of}
\left(
\mathop{\oplus}_{\ell=1}^{\infty}
S^{\ell}\right)
$$
where  $\ds  S^{\ell}$  is the $\ell$-th symmetric power of
$q^1,\ldots, q^n, p_1,\ldots, p_n$.  It holds $\ds 
\Pkt{S^{k}}{S^{\ell}} \subset S^{k+\ell-2}$, because the Poisson bracket 
satisfy
$ \Pkt{ q^i }{p_j} = - \Pkt{ p_j }{q^i} = \delta^{i}_{j} \quad
\text{others are 0}$, and   
$\ds  \hat{J} : \fraksp(2n,\mR)
\longrightarrow S^{2}$ is a Lie algebra isomorphism.

If we denote the dual of $\ds  S^{\ell}$ by 
$\ds  \kS{\ell}$, then we see the first 
cochain complex is    
\begin{align*}
        \ds  \cgf{1}{2n}{}{} & \cong \frakham_{2n}{}^{*} = 
        \text{the completion of }
\mathop{\oplus}_{\ell=1}^{\infty} \kS{\ell}\\
\noalign{and 
 the second
cochain complex is} 
\ds  \cgf{2}{2n}{}{}  & \cong \frakham_{2n}{}^{*} \wedge
\frakham_{2n}{}^{*} =
\mathop{\oplus}_{\ell=1}^{\infty} \kS{\ell}
\wedge 
\mathop{\oplus}_{k=1}^{\infty} \kS{k} 
 = 
\mathop{\oplus}_{1\le k\le \ell} \kS{k} \wedge  \frakS{\ell}
\\
\noalign{ and so on, now on we omit the comment of ``the completion'' for simplicity. }
\end{align*}
\label{section::weight}

\subsection{Weight of cochains} 

\begin{defn}[cf.\cite{KOT:MORITA}]
Define the \textbf{weight} of each element of
$\ds  \kS{\ell}$ by $\ell-2$.   

For each element of $\ds  
\kS{\ell_1} \wedge 
\kS{\ell_2} \wedge \cdots
\kS{\ell_s} $, define  its \textbf{weight} by $\ds  \sum_{i=1}^{s} (\ell_i
-2)$. 
\end{defn}

\begin{myRemark} 
Let $\ds  \sigma \in \kS{\ell} $ be a 1-cochain.  
Since $\ds  (\myd \sigma)(f_0,f_1) = -
\langle \sigma , \Pkt{f_0}{f_1}\rangle $,   
the contribution of $\sigma$ is when the case of  
$ \Pkt{f_0}{f_1}\in S^{\ell}$.  If 
$\ds  f_0 \in S^{p_0}$  and 
$\ds  f_1 \in S^{p_1}$, then  it must hold $p_0+p_1-2 = \ell$, namely, 
$\ds  \myd \kS{\ell} \subset \sum_{ p_0+p_1= 2+\ell}  \frakS{d_0}
\wedge \kS{d_1}$.  
Similarly, we see that 
$$ \myd (\kS{k}\wedge \frakS{\ell}) \subset \sum_{p_0+p_1+p_2=k+\ell+2}
\kS{p_0}\wedge  
\kS{p_1}\wedge  
\kS{p_2}.$$  
$p_0+p_1-2 = \ell$ is equivalent to 
$(p_0-2)+(p_1-2) = \ell-2$, and  $p_0+p_1+p_2=k+\ell+2$ is 
$
(p_0-2)+(p_1-2)+(p_2-2) =(k-2) +(\ell-2)
$.  
These show the reason of the definition of \textbf{weight} above. 
And 
we also see that the coboundary operator $\myd$ preserve the weight,
namely if a cochain $\sigma$ is of weight $w$, then 
$\ds  \myd (\sigma)$ is also of weight $w$.  
\end{myRemark}

Now we can decompose the cochain complex by
the weight 
$w$ as follows:  
$$ \cgf{\bullet}{2n}{}{w} := \sum_{w\text{-condition}} 
\Lambda ^{k_1} \kS{1} \otimes 
\Lambda ^{k_2} \kS{2} \otimes  \cdots $$
where $w$-condition is $ \ds  \sum_{j=1}^{\infty} (j-2) k_j= w$.  
$\ds   
\cgf{\bullet}{2n}{}{}  \cong \sum_{w=-2n}^{\infty}
\cgf{\bullet}{2n}{}{ w }
$.   
Since the coboundary operator $\myd$ preserves the weights and so we
have the natural splitting of cohomology groups like 
$$ 
\hgf{\bullet}{2n}{}{}  \cong \sum_{w=-2n}^{\infty}
\hgf{\bullet}{2n}{}{w }
$$
In order to investigate 
$\ds  
\hgf{m}{2n}{}{w}$, we have to handle  
$$ \cgf{m}{2n}{}{w} := \sum
\Lambda ^{k_1} \kS{1} \otimes 
\Lambda ^{k_2} \kS{2} \otimes  \cdots $$
$ \ds  \sum_{j=1}^{\infty}  k_j= m$ and   
$ \ds  \sum_{j=1}^{\infty} (j-2) k_j= w$. We have to be careful of no
contribution of ($\ds  k_2$)-term to $w$.     

Let 
$\ds  {\kham_{2n}^{0}}$ be the space of  
the Hamiltonian vector fields which vanish at the origin of $\ds \mR^{2n}$.  
Then we have 
$$
\ds  {\kham_{2n}} \cong  \mathop{\oplus}_{\ell=1}^{\infty}S^{\ell},\quad   
\ds  {\kham_{2n}^{0}} \cong 
\mathop{\oplus}_{\ell=2}^{\infty}S^{\ell}
\quad\text{and}\quad  
\ds  {\kham_{2n}^{1}} \cong 
\mathop{\oplus}_{\ell=3}^{\infty}S^{\ell} 
\ . $$   
We look for 
the relative 
$\ds  
\HGF{\bullet}{2n}{0}{w} $.

\begin{myRemark}[\cite{KOT:MORITA}]
If $w$ is odd then 
$\ds  
\CGF{\bullet}{2n}{}{w}  = \{ 0\}$. 
\end{myRemark}
$\ds  
\begin{pmatrix} -I & 0 \\ 0 & -I \end{pmatrix}
=: - \; \text{Id}  
$ is an
element of $Sp(2n,\mR)$ and 
each cochain $\sigma$ should be invariant under the action of
$-\text{Id}$.  
$$\forall \sigma \in \text{C}_{GF}^{m}({\kham_{2n}}) _ w := \sum
\Lambda ^{k_1} \kS{1} \otimes 
\Lambda ^{k_2} \kS{2} \otimes  \cdots \quad\text{, where}\quad 
 \sum_{j=1}^{\infty}  k_j= m \quad\text{and}\quad    
\sum_{j=1}^{\infty} (j-2) k_j= w$$ 
we see that
\begin{align*}
\sigma = &  (-\text{Id}) \cdot \sigma 
=  (-1)^{k_1} (+1)^{k_2} (-1)^{k_3} (+1)^{k_4} \cdots \sigma \\
= & (-1)^{ k_1 +  k_3 + k_5 +  k_7 + \cdots }  \sigma 
=  (-1)^{ -k_1 + k_3 + 3 k_5 + 5 k_7 +  \cdots }  \sigma \\
= & (-1)^{ -k_1 + 0 k_2 + k_3 + 
+2 k_4 3 k_5 + 4k_6 + 5 k_7 +  \cdots }  \sigma 
= (-1) ^{w}  \sigma 
\end{align*}
thus $\sigma =0$ if $w$ is odd.  
\kmqed

\begin{myRemark} When $n=1$, 
the {\em type} of Metoki
(\cite{metoki:shinya}) 
and the {\em weight} here are related by 2 times of {\em type} is equal to 
{\em weight}.
\end{myRemark}
We remember that 
\begin{align*}
        \cgf{m}{2n}{0}{w} & = \sum
\Lambda ^{k_1} \kS{1} \otimes 
\Lambda ^{k_2} \kS{2} \otimes  \cdots \quad\text{, where}\  
 \sum_{j=1}^{\infty}  k_j= m  ,\   
\sum_{j=1}^{\infty} (j-2) k_j= w\ \text{and}\ k_1 = 0      
\\
        \cgf{m}{2n}{1}{w} & = \sum
\Lambda ^{k_1} \kS{1} \otimes 
\Lambda ^{k_2} \kS{2} \otimes  \cdots \quad\text{, where}\  
 \sum_{j=1}^{\infty}  k_j= m ,\  
\sum_{j=1}^{\infty} (j-2) k_j= w\ \text{and}\ k_1 = k_2 = 0     
\end{align*}

\subsection{Relativity}
On our original Lie algebra of Hamiltonian vector fields $\ds  
\Hd_f$ of polynomial potential functions $f$, we have the natural group
action of $\ds  Sp(2n,\mR)$:
$$ a\cdot \Hd_{f} = \Hd_{f\circ a^{-1}},\qquad a\cdot f := f \circ
a^{-1} \ . $$

\begin{defn} On the exterior algebra of the space of polynomial
functions of $\ds  \mR^{2n}$, the group 
$Sp(2n,\mR)$ acts naturally by 
$$ a\cdot (f_1 \wedge f_2 \wedge \cdots f_m) := 
(a \cdot f_1) \wedge \cdots \wedge (a \cdot f_m)$$
where $\ds  f_j$ are polynomials. 
We define the infinitesimal action of $\xi\in
\fraksp(2n,\mR)$ by   
$$ \xi \cdot  (f_1 \wedge f_2 \wedge \cdots f_m) 
:= \frac{d}{dt} \left( \exp(t\xi)\cdot  (f_1 \wedge f_2 \wedge \cdots
f_m)\right) 
_{\mid t=0}  \ . 
$$
\end{defn}
\begin{myRemark}
$\ds  f_1 \wedge f_2 \wedge \cdots f_m $ 
looks strange, but the corresponding real object is 
$\ds  \Hd_{-f_1} \wedge \Hd_{-f_2} \wedge \cdots \wedge \Hd_{-f_m}$. 
Thus the degree of $f_j$ is 1, and so 
$\ds  f \wedge g = - g \wedge f$ for each functions $f,g$.  
\end{myRemark}

\begin{kmProp}
The above infinitesimal action by $\ds  \xi\in\fraksp(2n,\mR)$
is a derivation of degree 0 and 
$$ \xi \cdot f = \Pkt{ \hat{J}(\xi)}{f} \quad \text{for each
polynomial}\ f \ .$$ 
\end{kmProp}
Proof: 
The first assertion of being derivation of degree 0 is trivial. The
second assertion is:

$$ \xi \cdot f := \frac{d}{dt} \exp(t\xi)\cdot f _{\mid t=0}  
= \frac{d}{dt} f\circ \exp(-t\xi) _{\mid t=0} = 
\langle df, - \xi_M \rangle = \omega( \Hd_f, -\Hd_{ \hat{J}(\xi)}) =
\Pkt{\hat{J}(\xi)}{f} . $$
This corresponds to 
$\ds  \xi \cdot \Hd _{f} := 
\frac{d}{dt} \exp(t\xi)\cdot \Hd_{f} {}_{\mid t=0} = {\cL}_{-\xi_{M}} \Hd_{f}
= -[ \Hd_{ \hat{J}(\xi)}, \Hd_{f}] = \Hd_{ \Pkt{ \hat{J}(\xi)}{f}}$. 
\kmqed

\bigskip

To determine relative cochain complex 
$\ds  
\CGF{\bullet}{2n}{0}{}$,  
there are two  conditions, one is 
$\ds  i_{\xi} \sigma =0$, and the other is 
$\ds  i_{\xi} \myd \sigma =0$  
for each  $m$-cochain $\ds  \sigma $,    
where 
$\ds  i_{\xi}$ is 
the interior product with respect to 
$\xi \in \frakk = \fraksp(2n,\mR)= S^2$.  
Since $\ds  i_{\xi}$
is a 
skew-derivation of   
degree $-1$,  
in order to know the effect of 
$\ds  i_{\xi}$, it is enough to 
know the operation of $\ds  
i_{\xi} \sigma$ for 1-cochain $\sigma$.  Going back to Hamiltonian
vector fields, we see that 
$$ i_{\xi}\sigma = \langle \sigma , \xi_{M}\rangle = 
 \langle \sigma , \Hd_{ \hat{J}(\xi)} \rangle = 
 \langle \sigma , -  \hat{J}(\xi) \rangle . $$ 
In general, it holds 
$$ 
\ds  ( i_{\xi}\sigma)(f_1,f_2,\ldots) = \sigma(-
\hat{J}(\xi),f_1,f_2, \ldots) \quad \text{for}\quad \sigma \in
\cgf{m}{2n}{}{}  . $$

\begin{kmProp}\label{int:xi}
Let $\ds  
\sigma \ne 0$ and   
$\ds  
\sigma \in 
\Lambda^{k_1} \kS{1} \otimes 
\Lambda^{k_2} \kS{2} \otimes
\Lambda^{k_3} \kS{3} \otimes
\cdots $ (where $\ds  \sum
k_j =m$). If $\ds  i_{\xi} \sigma = 0$ for $\forall \xi \in
\frakk = \fraksp(2n,\mR)$, then 
$\ds  
\sigma \in 
\Lambda^{k_1} \kS{1} \otimes 
\Lambda^{k_3} \kS{3} \otimes 
\Lambda^{k_4} \kS{4} \otimes \cdots $, i.e., $\ds k_2=0$.  
\end{kmProp}

\medskip

Proof:  Since $\ds  i_{\xi} \sigma = - \la \sigma,
\hat{J}(\xi)\ra $ for 1-cochain and  
$\ds  
\hat{J}(\xi) \in S^{2} $,   
we see that 
$\ds  i_{\xi} \sigma =0$ if 
$\ds  \sigma  \in \kS{\ell} $ with $\ell \ne 2$.  
If 
$\ds  \sigma  \in \kS{2} $ satisfies  
$\ds  i_{\xi} \sigma =0$ for $\forall \xi \in\fraksp(2n,\mR)$, 
then  
$\ds  \sigma =0 \in \kS{2} $ because 
$\ds  S^{2} = \fraksp(2n,\mR)$ is semi-simple and 
$\ds  \hat{J}(\xi)$ generate $\ds  S^{2}$.   
Now we may rewrite 
$$ \sigma = \sum_{A} \tau_A \wedge \rho_A$$ where 
$\ds  \tau _A \in 
\Lambda^{k_2} \kS{2}$,  
$\ds  \rho_A \in 
\Lambda^{k_3} \kS{3} \otimes 
\Lambda^{k_4} \kS{4} \otimes \cdots $ and 
$\ds  \rho_A $ are linearly independent. 
Using the fact $\ds  i_{\xi} \rho_A =0$ 
for $\forall \xi\in\frakk$, we have 
$\ds  0= i_{\xi} \sigma = \sum_{A} 
 i_{\xi} \tau_A \otimes  \rho_{A} 
$ and those imply that 
$\ds  i_{\xi} \tau_A = 0 $ 
for $\forall \xi\in\fraksp(2n,\mR)$ and $A$. If $0< k_2 \le \dim \frakS{2}=
n(2n+1)$ then $\tau_A=0$ for $\forall A$ if necessary we can use 
$\ds   i_{\xi_{k_2}} \cdots  i_{\xi_1}
 \tau_A = 0 $. Thus $ \sigma = 0$.  
\kmqed

\bigskip

%%A
Thus, 
the condition 
$\ds  i_{\xi}\sigma=0$ for $\xi\in \fraksp(2n,\mR)$ 
imply that $\ds  \kS{2}$ does not appear
in 
$\ds  \CGF{\bullet}{2n}{0}{}$. 

\kmcomment{
and 
\begin{align} 
\CGF{m} {2n}{0}{} 
& = \Lambda^{m} (
\kS{ 3} \oplus \cdots)^{Sp(2n,\mR)} 
= \{ \sigma 
         \in \cgf{m}{2n}{1}{}  \mid 
 Sp(2n,\mR)\text{-invariant} \} 
\\
\noalign{and so}
\CGF{m}{2n}{0}{w}
&= 
\{ \sigma 
         \in \cgf{m}{2n}{1}{w}  \mid 
 Sp(2n,\mR)\text{-invariant} \} \label{sub:Kont}
\end{align}
}
%BB

%\begin{myRemark}
The other condition of being relative cochain is 
$\ds  i_{\xi} \myd \sigma =0$  
for each $m$-cochain $\ds  \sigma $.   Again going back to
Hamiltonian vector fields, for each 1-cochain $\sigma$, we see that 
\begin{align*} 
\langle i_{\xi}\myd \sigma, \Hd_f \rangle  =& 
( \myd \sigma)  (\xi_{M}, \Hd_f)  = - 
 \langle \sigma , [ \xi_{M}, \Hd_f] \rangle = 
 - \langle \sigma , [ \Hd_{ \hat{J}(\xi)},\Hd_f] \rangle = 
 \langle \sigma ,  \Hd_{\Pkt{ \hat{J}(\xi)}{f}} \rangle 
 \\
\noalign{and so}  
\langle i_{\xi}\myd \sigma, f \rangle  =& 
 \langle \sigma ,  \Pkt{ \hat{J}(\xi)}{f} \rangle 
=  
 \langle \sigma ,  \xi \cdot f  \rangle 
%= 
% \langle - \xi \cdot \sigma ,   f  \rangle 
\end{align*}
for
 1-cochain $\sigma$, and $ \xi \in \fraksp(2n,\mR),
 f\in \text{C}^{\infty}(M)$.

For m-cochain 
$\ds  \ds  \tau = 
\sigma_1 \wedge \cdots \wedge
\sigma_m$ ($\sigma_j$ are 1-cochains), 
since $\myd$ is a skew-derivation of degree $+1$, we have  
$$ 
\myd \tau   =  \sum_{j=1}^m (-1)^{j+1} 
\sigma_1 \wedge \cdots \wedge \myd \sigma_j \wedge \cdots 
\wedge \sigma_m  , 
$$
and since the interior product $\ds  i_{\xi}$ for 
each $\xi \in \fraksp(2n,\mR)$ is a skew-derivation of degree
$-1$, we have   
\begin{align*}
\left( i_{\xi} \circ\myd \right)\tau  =& i_{\xi} \sum_{j} (-1)^{j+1} 
\sigma_1 \wedge \cdots \wedge \myd \sigma_j \wedge \cdots 
\wedge \sigma_m \\
=& 
\sum_{\ell < j} (-1)^{j+1} (-1)^{\ell+1} 
\sigma_1 \wedge \cdots \wedge 
i_{\xi} \sigma_\ell \wedge   
\cdots \wedge \myd \sigma_j \wedge \cdots 
\wedge \sigma_m \\
& + 
\sum_{j} (-1)^{j+1} (-1)^{j+1} 
\sigma_1 \wedge 
\cdots \wedge \left((i_{\xi}\circ \myd )\sigma_j\right) \wedge \cdots 
\wedge \sigma_m 
\\
& + 
\sum_{\ell > j} (-1)^{j+1} (-1)^{\ell+2} 
\sigma_1 \wedge \cdots 
\wedge \myd \sigma_j 
\wedge \cdots 
\wedge i_{\xi} \sigma_\ell \wedge  \cdots 
\wedge \sigma_m 
\end{align*}
If $\ds  i_{\xi} \sigma_j=0$ ($j=1,\ldots,m$), then we have 
$$ 
\left( i_{\xi} \circ\myd \right)\tau   
= 
\sum_{j} \sigma_1 \wedge 
\cdots \wedge \left( (i_{\xi} \circ \myd) \sigma_j\right) \wedge \cdots 
\wedge \sigma_m ,
$$ namely, 
$\ds  
 i_{\xi} \circ\myd $ is a derivation of degree 0.    
Thus, under the condition of $\ds  i_{\xi} \sigma =0$ for any
cochain $\sigma$ and $\ds 
\xi\in  \fraksp(2n,\mR)$, 
$\ds  i_{\xi}\circ \myd $ becomes  
an ordinary derivation of degree 0.  

\begin{kmProp}
Let $\sigma$ be a m-cochain with 
$\ds  i_{\xi} \sigma  = 0$ for $\forall
\xi\in\fraksp(2n,\mR)$. Then 
$\ds  (i_{\xi}\circ \myd) 
  $  behaves as derivation of degree 0 and 
characterized by 
$\ds  \la (i_{\xi}\circ \myd) \sigma, f\ra  =
\la  \sigma , \Pkt{\hat{J}(\xi)}{f}\ra $ 
for each 1-cochain $\sigma$ and $\xi\in \fraksp(2n,\mR)$. 
\end{kmProp} 

\begin{myRemark}
It may be a better way to recall Cartan's formula $\ds  {\cL}_{\xi }
= \myd \circ i_{\xi} + i_{\xi}\circ \myd$ in order to prove the above.    
\end{myRemark}

%CC
From Proposition \ref{int:xi} and the discussion above, we see that 
%%A
\begin{align} 
\CGF{m} {2n}{0}{} 
& = \Lambda^{m} (
\kS{ 3} \oplus \cdots)^{Sp(2n,\mR)} 
= \{ \sigma 
         \in \cgf{m}{2n}{1}{}  \mid 
 Sp(2n,\mR)\text{-invariant} \} 
\\
\noalign{and so}
\CGF{m}{2n}{0}{w}
&= 
\{ \sigma 
         \in \cgf{m}{2n}{1}{w}  \mid 
 Sp(2n,\mR)\text{-invariant} \} \label{sub:Kont}
\end{align}
%BB

%HERE2
The group $Sp(2n,\mR)$ acts on cochain complexes as the dual action 
of that  
of on the exterior algebra of polynomial functions on
$\ds  \mR^{2n}$. Thus, the precise definition is 
\begin{align*} 
\left( a\cdot \sigma \right) (f_1,\ldots, f_m)  &= 
\sigma \left(   a^{-1} \cdot f_1, \ldots , 
  a^{-1} \cdot f_m \right) 
= 
\sigma \left(   f_1 \circ a_M , \ldots , f_m \circ a_M\right) \\ 
\la a\cdot \sigma , f_1 \wedge \ldots \wedge f_m \ra &= 
\la  \sigma , a^{-1} \cdot ( f_1 \wedge \ldots \wedge f_m ) \ra 
= 
\la \sigma,    f_1 \circ a_M  \wedge \cdots \wedge  f_m \circ a_M\ra   
\end{align*}
for a general m-cochain $\sigma$.  

\begin{myRemark} We re-confirm here that our coboundary operator $\myd$
is compatible with the group-action of $\ds  Sp(2n, \mR) $,
i.e,.    $\ds  a \cdot \myd = \myd \circ a_M $ holds 
for each symplectic automorphism $a$. 
It is enough only to show for 1-cochain $\sigma$.    
We see that 
\begin{align*} \ds  
\left(\myd (a \cdot \sigma)\right)(f,g) &= 
- \langle a \cdot \sigma , \Pkt{f}{g}\rangle 
= - \langle  \sigma , \Pkt{f}{g} \circ a^{-1}_{M} \rangle 
= - \langle  \sigma , \Pkt{f\circ a^{-1}_M }{g\circ a^{-1}_{M}} \rangle \\&
= ( \myd \sigma ) ( f\circ a^{-1}_M ,g\circ a^{-1}_M) 
=  ( \myd \sigma )  ( a\cdot f , a\cdot g)   .
\end{align*}
for $\ds  \forall 
a\in Sp(2n,\mR) $.   
\end{myRemark}

\begin{defn}
We define the infinitesimal action of $\ds 
\xi\in\fraksp(2n,\mR)$ for each cochain $\sigma$ by 
$$ \ds  \xi \cdot \sigma = \frac{d}{dt} \exp(t\xi)\cdot \sigma
_{\mid t=0} \ .$$ 
\end{defn}
\begin{kmProp}\label{prop::advantage}
The infinitesimal action $\xi\in\fraksp(2n,\mR)$ preserves the cochain
complex $ 
\ds  \cgf{m}{2n}{0}{} 
$, behaves as an ordinary derivation of degree 0 and 
$$ \la \xi \cdot \sigma , f \ra = - \la \sigma , \xi \cdot f\ra =
- \la \sigma, \Pkt{ \hat{J}(\xi)}{f}\ra = - \la ( i_{\xi}\myd \sigma,
f\ra \ .$$
Thus, $\ds   (i_{\xi}\myd)(\sigma) = - \xi\cdot \sigma$ holds
for each m-cochain $\sigma$ with 
$\ds   i_{\xi}\sigma = 0$ ($\xi\in\fraksp(2n,\mR)$). 
\end{kmProp}

An advantage of the relation in Proposition \ref{prop::advantage} is
that $\ds  \xi :
\cgf{m}{2n}{0}{} \longrightarrow 
\cgf{m}{2n}{0}{} $ is a derivation of degree 0 with respect to
the wedge product but also a derivation of degree 0 inside of 1-cochain. 
Namely, each 1-cochain is a linear combination of symmetric powers and 
since $\ds   f \mapsto 
\Pkt{ \hat{J}(\xi)}{f}$ is a derivation for each polynomial function
$f$, we have 
\begin{align*} \xi \cdot 
\left( e_1{}^{k_1} \cdots 
e_{\ell}{}^{k_{\ell}} \cdots 
e_{2n}{}^{k_{2n}}\right)  
& = \sum_{\ell}
 e_1{}^{k_1} \cdots 
\xi \cdot \left(e_{\ell}{}^{k_{\ell}}\right) \cdots 
e_{2n}{}^{k_{2n}}
\\
&
= \sum_{\ell}
 e_1{}^{k_1} \cdots 
\cdot \left( k_{\ell} e_{\ell}{}^{k_{\ell}-1} \xi( e_{\ell}) \right) \cdots 
e_{2n}{}^{k_{2n}} \ .
\end{align*}

\section{Decomposition of $m$-th cochain complex of weight $w$} 
Our concern in this section is, given a pair of positive integers $(m,w)$, find all
possibilities of cochain complex by denoting the sequences $(k_3,k_4,\ldots)$
of non-negative integers of multiplicity which satisfies  
$\ds  \sum_{j\ge 3} k_j = m$ and 
$\ds  \sum _{j\ge 3} (j-2)k_j=w$. 
We have to be careful about 
$\ds  \Lambda^{k_{\ell}} \kS{\ell} = \{ 0\}$ may be happen
when $\ds  k_{ell} > \dim \kS{\ell}$, where 
$\ds  \dim \kS{\ell} = (\ell + 2n-1)!/( \ell! (2n-1)!) $.    

\renewcommand{\hk}{\hat{k}} %\newcommand'Å'Í'È'­\renewcommand'ɕύX

By shifting the
indices by $-2$, we rearrange our situation as below.  
Given a pair of non-negative integers $(m,w)$, we would like to find all
sequences $(\hk_1,\hk_2,\ldots)$ of non-negative integers satisfying 
\begin{align}
%(1) 
\ds  \sum_{j\ge 1} \hk_j = m \label{A:eqn} \quad \text{and} \quad 
%and (2) 
\ds  
\sum _{j\ge 1} j \hk_j=w
%\label{B:eqn} 
\end{align} 
$m \le w$ is a necessary condition and denote the set of sequences
$(\hk_1,\hk_2,\ldots)$ satisfying (\ref{A:eqn}) by $F(m,w)$, 
By subtracting the first equation of (\ref{A:eqn}) from
the second one of (\ref{A:eqn}), we have 
$$
\sum _{j\ge 2} (j-1) \hk_j=w-m \ . 
$$ Thus, we have 
the all solutions of (\ref{A:eqn}),   
$ F(m,w) $  
by 
the next recursive formula: 
$$ F(m,w) = \mathop{\bigsqcup}_{\hk_1=\max(0,2m-w)}^{m} \{ (\hk_1, x) \mid x 
\in F(m-\hk_1, w-m)\}
$$

%NG in quote-env  \verb|gkf_act-1.mpl|
%OK \verbatiminput{gkf_act-1.mpl}
For a given $(m,w)$, we have  
the Maple script \texttt{gkf\_act-1.mpl}, which shows us    
the all solutions of (\ref{A:eqn}),   
%KM \input{gkf_20.table}

There is a primitive question ``do the solutions exist
for all $m \le w$? or how many?''. 
We will give an answer to this question, here. We join our two equations
into one equation as below:

\begin{align*}
w &= \hk_1+ 2 \hk_2+ \cdots + s \hk_s \\
  &= (\underbrace{1+\cdots+1}_{\hk_1})
  + (\underbrace{2+\cdots+2}_{\hk_2}) + \cdots +
   (\underbrace{s+\cdots+s}_{\hk_s}) \\
&= 
   (\underbrace{s+\cdots+s}_{\hk_s}) + \cdots  
  + (\underbrace{2+\cdots+2}_{\hk_2}) + 
  (\underbrace{1+\cdots+1}_{\hk_1})\\
&= \ell_1+ \ell_2+ \cdots + \ell_m 
\end{align*}
where $\ell_1 \ge \ell_2 \ge \cdots \ge \ell_{m} \ge 1$.  
This is a partition of $w$ with length $m$ or a Young
diagram of height $m$  
with $w$ cells.   
Conversely, for a partition of $w$ 
\begin{equation}\label{eqn:B}
\begin{split}
& w = \ell_1+\ell_2+\cdots+\ell_m \\
& \ell_1 \geqq \ell_2 \geqq \cdots \geqq \ell_m \geqq  1
\end{split}
\end{equation}
$ \hk_i := \# \{ j \,|\, \ell_j=i \}$ gives   
a solution of (\ref{A:eqn}).   
That means there is a one-to-one correspondence between the  
solution of (\ref{A:eqn}) and all the partitions of $w$ with length $m$ or a Young
diagram of height $m$  
with $w$ cells.  
\begin{myRemark}
Be careful of difference of definition of 
$\ds   \hk_i := \# \{ j \,|\, \ell_j=i \}$ 
and the conjugate Young diagram 
$\ds  c_i :=  \# \{ j \,|\, \ell_j \geqq i \}$.  
\end{myRemark}

\begin{kmProp}

By $r(m,w)$ we mean the number of solutions of (\ref{A:eqn}).  
\bigskip

For integers $m>0$, $w\geqq 0$, we define $\tilde{r}(m,w)$ the number of
solutions of 
\begin{equation}\label{eqn:C}
\begin{split}
& w = \ell_1+\ell_2+\cdots+\ell_m \\
& \ell_1 \geqq \ell_2 \geqq \cdots \geqq \ell_m \geqq  0
\end{split}
\end{equation}
By an elementary observation below 
$$ 
w  = \ell_1+ \ell_2 + \cdots + \ell_m \quad \text{with} \quad \ell_i \ge 0 \quad
\text{and} 
\quad  w+m  = (\ell_1+1)+ (\ell_2+1) + \cdots + (\ell_m+1) \quad \text{with}\quad
 (\ell_i+1) > 0 $$
we see $\tilde{r}(m,w)=r(m,w+m)$.  
The generating function of $\tilde{r}(m,w)$ is known as 
$$ \sum_{c=0}^\infty \tilde{r}(m,c)x^c = \prod_{k=1}^m \left( \frac{1}{1-x^k} \right).$$
\end{kmProp}

%HERE
\begin{exam}\label{exam::a}
Let us try to know all possibilities when weight=2 case.  

When $m=1$, then $\ell_1 = 2$, so we have $\hk_2= 1$ and $\hk_j=0$ ($j
\ne 2$). Thus, $k_4=1$.   

When $m=2$, then $ 2 = \ell_1 + \ell_2$ ( $\ell_1 \ge \ell_2 \ge 1$),  
$\ell_1 = \ell_2 =1$, so we have $\hk_1 = 2$ and 
$\hk_j=0$ ($j \ne 1$). Thus, $k_3 = 2$. These show  
\begin{alignat*}{3}
& \cgf{1}{2n}{1}{2} 
= \kS{4} 
& \CGF{1} {2n}{0}{2} 
&= \kS{4}^{Sp(2n,\mR)} (= \{\mathbf{0} \})\\
& \cgf{2} {2n}{1}{2} 
= \Lambda^2 \kS{3} \qquad 
& \CGF{2}{2n}{0}{2} 
&= (\Lambda^2 \kS{3})^{Sp(2n,\mR)}  
\end{alignat*} 
In the above, $\ds  \kS{4}$ is an irreducible representation
of $Sp(2n,\mR)$ and so 
$\ds  \kS{4}^{Sp(2n,\mR)} = \{ \mathbf{0} \}$. Also,
concerning 
$\ds  \Lambda^2 \kS{3}$, if n=1 we know  
$\ds  \Lambda^2 \kS{3} = \frakS{0} \oplus \frakS{4}$ as we
will see in Example \ref{exam::d}. If 
n=2 then by the help of Littlewood-Richardson rule, 
we get a little complex expression 
$\ds  \Lambda^2 \kS{3} = \frakS{0} \oplus \frakS{4}
\oplus V_{<(1,1)>} 
\oplus V_{<(2,2)>} 
\oplus V_{<(3,3)>} 
\oplus V_{<(5,1)>} 
$, where $\ds  V_{<(p,q)>}$ is the irreducible representation of 
the natural action of $Sp(4,\mR)$ on $\ds  \mR^4$,
corresponding the Young diagram $(p,q)$, and $\ds  V_{<(p,0)>} =
\kS{p}$.  Thus, if n=1 or 2, then we see  
$\ds    
\dim \CGF{2}{2n}{0}{2} =1$.
Actually, we get a lot of help from representation theory in this
project.  
\end{exam}

\begin{exam} weight=4 case:  \label{exam::b}

When $m=2$, i.e., $ 4 = \ell_1 + \ell_2$ ( $\ell_1 \ge \ell_2 \ge 1$),
then   
$(\ell_1,\ell_2) = (3,1) $ or  $(2,2)$, 
so we have $(\hk_1 = 1, \hk_3 = 1) $, or   
$\hk_2=2$ . Thus, $(k_3 = 1, k_5=1)$ or $(k_4=2)$. 
When $m=3$, i.e., $ 4 = \ell_1 + \ell_2 +\ell_3$ ( $\ell_1 \ge \ell_2
\ge \ell_3 \ge  1$), $\ell_1=2, \ell_2 = 1, \ell_3=1$. Thus 
$(\hk_2=1, \hk_1=2)$, so $(k_3=2, k_4=1)$.  

These show  
\begin{alignat*}{2}
& \cgf{1} {2n}{1}{4} = \kS{6} 
& \cgf{2} {2n}{1}{4} &= 
( \kS{3} \otimes \frakS{5} )\oplus 
\Lambda^2 \kS{4}
\\ 
& \cgf{3}{2n}{1}{4} 
= \Lambda^2 \kS{3} \otimes \kS{4}   \qquad 
& \cgf{4}{2n}{1}{4} & = \Lambda^4 \kS{3}
\end{alignat*}
\end{exam}

\begin{exam} \label{exam::c}
By the same way, we have weight=6 case: 
\begin{align*} 
\cgf{1}{2n}{1}{6} & = \frakS{8}  
\\
\cgf{2}{2n}{1}{6} 
& =  
\left(  \kS{3} \otimes \frakS{7}  \right)
\oplus \left(  \kS{4} \otimes \frakS{6}  \right)
\oplus \Lambda^2 \kS{5}  
\\ 
\cgf{3}{2n}{1}{6} 
& = \left( \Lambda^2 \kS{3} \otimes \frakS{6}  \right)
\oplus \left(  \kS{3} \otimes \frakS{4} \otimes \frakS{5}  \right)
\oplus  \Lambda^{3} \kS{4}  
\\ 
\cgf{4}{2n}{1}{6} 
& = \left( \Lambda^3 \kS{3} \otimes \frakS{5}  \right)
\oplus \left( \Lambda^{2}  \kS{3} \otimes \Lambda^{2} \frakS{4} \right)
\\
\cgf{5} {2n}{1}{6} 
& =  \Lambda^4 \kS{3} \otimes \frakS{4}  
\\
\cgf{6} {2n}{1}{6} 
& = \Lambda^{6} \kS{3} 
\end{align*}  
If $n=1$, then $\dim \kS{3} = 4$ and so we have 
$ \cgf{6} {2}{1}{6} = \{ \mathbf{0} \} $.  
\end{exam}

We express the data which we got by the next table in short form, dividing
into direct sum components. Our abbreviation
rule is that 
\begin{enumerate}
\item only pick up the $i$'s with $k_i>0 $, 
\item 
if $k_i>1 $ then express the multiplicity by the power like 
$\ds  i^{k_i}$, i.e.,  
$\ds  \left( i^{k_i}\ j^{k_j}\ \cdots \right)$  
\item 
if $k_i=1$ then  
only write 
$\ds  i$ 
\end{enumerate}
Using the rule, 
Example \ref{exam::c} above with weight=6 can be written in the next
table: 

\hfil
\begin{tabular}[t]{|c|c|l|}
\hline
degree & ref.\# & type \\\hline
1 & 1 & ($8$) \\\hline
2 & 1 & ($3\ 7$)\\
  & 2 & ($4\ 6$)\\
  & 3 & ($5^2 $)\\\hline
\end{tabular}
\hfil 
\begin{tabular}[t]{|c|c|l|}
\hline
degree &ref.\# & type \\\hline 
3 & 1 & ($3^2\ 6$)\\
  & 2 & ($3\ 4\ 5$)\\
  & 3 & ($4^3 $)\\\hline
\end{tabular}
\hfil  
\begin{tabular}[t]{|c|c|l|}
\hline
degree &ref.\# & type \\\hline
4 & 1 & ($3^3\ 5$)\\
  & 2 & ($3^2\ 4^2$)\\\hline
5 & 1 & ($3^4\ 4$)\\\hline
6 & 1 & ($ 3^6$) \\\hline
\end{tabular}
\hfil 

\bigskip

\section{Aid from Representation theory of $Sp(2n,\mR)$} 
In order to compute the relative cohomology groups, we have to know some basis
of the cochain complex and the concrete matrix representation of the
coboundary operator, and that rank.  Fortunately, in the case of
$\ds \CGF{m} {2n}{0}{w} $,
there is a very sophisticated way to know the dimension without knowing basis.  
Namely, wan consider the natural representation of $\ds 
Sp(2n,\mR)$ acting on $\ds  \mR^{2n}$.  Then the all
irreducible representations are known by Young diagram of length at most
$n$. Also, the $p$-th symmetric tensor product of $\mR^{2n}$ is an
irreducible representation and is identified with $\ds 
\kS{p}$ or its dual in this paper.    
The corresponding Young diagram is $(p) =
\underbrace{
\fbox{\phantom{x}}\, 
\fbox{\phantom{x}}\, 
\cdots \, 
\fbox{\phantom{x}}\, 
\fbox{\phantom{x}}}_{p-\text{times}}$.

Since  we have (\ref{sub:Kont}), 
if we know 
the index of the trivial representation in 
$\ds  
\cgf{m}{2n}{1}{w} 
$ 
by some help of representation theory of $Sp(2n,\mR)$, 
that is just equal to 
$\ds  
\dim \CGF{m}{2n}{0}{w} 
$.  Thus, we can calculate the Euler characteristic number without
knowing cohomology data at least theoretically.  
When $n=1$, the representation theory is rather clear and that is very
helpful as we see later.  

\subsection{Assistance of $Sp(2,\mR)$} 
Instead of the variables $\ds  (q^1,p_1)$ in $\mR^2$, 
we use the classical notation
$(x,y)$, and let $ \tz{a}{k}$ be $\ds  \dfrac{x^{a}}{a!}
\dfrac{y^{k-a}}{(k-a)!} \in S^{k}$ 
and the dual basis in $\ds  \kS{k}$  of $\ds 
\tz {a}{k} $ by the notation $\nz{a}{k}$.   
Now, the Poisson bracket $\ds  \Pkt{f}{g} $ is 
the Jacobian $\ds \frac{\pdel (f,g)}{\pdel (x,y)}$ and get 
the relations 
$\ds  \Pkt{xy}{x^2}= -2 x^2$, 
$\ds  \Pkt{xy}{y^2}=  2 y^2$, 
$\ds  \Pkt{x^2}{y^2}= 4 x y$ 
and so we have the correspondence below with the famous matrices through
the momentum mapping $J$ of the natural symplectic action $Sp(2,\mR)$: 
$$ 
 xy \leftrightarrow H = \begin{pmatrix} 1 & 0 \\0 & -1\end{pmatrix}, 
\frac{y^2}{2} \leftrightarrow X = \begin{pmatrix} 0 & 1 \\0 & 0\end{pmatrix}, 
-\frac{x^2}{2} \leftrightarrow Y = \begin{pmatrix} 0 & 0 \\1 & 0\end{pmatrix} 
$$
and we have already seen that 
the component 
$\ds  S^2$ is a subalgebra of 
$\ds  {\kham_{2}^{0}}$ and isomorphic to 
$\ds  {\fraksp}(\mR^2, \omega) \cong {\fraksl}(2,\mR)$.

It is well-known that $\ds  S^{\ell}$ or the (symplectic) dual space
$\ds  \kS{\ell}$ are the irreducible representations of
$Sp(2,\mR)$ and Clebsch-Gordan rule, for instance, an irreducible
decomposition of tensor product is the following:   
$$ \kS{k}\otimes \frakS{\ell} = \frakS{k+\ell} \oplus \frakS{k+\ell-2} \oplus \cdots \oplus \frakS{|k-\ell |}$$ 

Since $\ds  
\text{C}_{GF}^{m} (\kham_{2}^{0}  ) _w 
= \sum_{\sum k_j = m,\  \sum(j-2)k_j=w}
\left(\Lambda^{k_2} \kS{2} \otimes 
\Lambda^{k_4} \kS{4} \otimes \cdots \right)
$, if  
we know the irreducible decomposition of 
$\ds  \Lambda^{k}\kS{\ell}$, then after tensor product-ing
those, we 
have the complete irreducible decomposition, 
and can pick up the trivial representation,    

\bigskip

\begin{exam} \label{exam::d}
As shown in Example \ref{exam::a}, $\ds  \Lambda^2 \kS{3}$
is a component of 
$\ds  
\cgf{2}{2n}{1}{2}$, we decompose  
$\ds  \Lambda^2 \kS{3}$ into irreducible components when
$n=1$.  
$\nz{\ell_1}{3}  \wedge 
\nz{\ell_2}{3} $ ($0\le \ell_1 <  \ell_2\le 3$) 
are a basis of $\ds  \Lambda ^{2} \kS{3}$.  
We will find the weight vector space from 
$$ T = \sum_{ 0 \le \ell_1 < \ell_2 < 3} 
c_{ \ell_1, \ell_2 } 
\nz{\ell_1}{ 3}  \wedge 
\nz{\ell_2}{ 3}
$$ 
which must be zero by the dual action of $\ds  X$.  
Since our action is 
$\ds  \nz{0}{3}\mapsto  3\nz{1}{3}$, 
$\ds  \nz{1}{3}\mapsto  2 \nz{2}{3}$, 
$\ds  \nz{2}{3}\mapsto  1 \nz{3}{3}$, and 
$\ds  \nz{3}{3}\mapsto  0$, we have 
\begin{align*}  
0 = &  c_{01}(3 \nz{1}{3} \wedge \nz{1}{3} + \nz{0}{3} \wedge 2 \nz{2}{3} ) 
+ c_{02} ( 3\nz{1}{3} \wedge \nz{2}{3} + \nz{0}{3} \wedge  \nz{3}{3} ) 
+ c_{03} ( 3\nz{1}{3} \wedge \nz{3}{3} + 0 ) \\& 
+ c_{12} ( 2 \nz{2}{3} \wedge \nz{2}{3} + \nz{1}{3} \wedge  \nz{3}{3} ) 
+ c_{13} ( 2 \nz{2}{3} \wedge \nz{3}{3} + 0 )
+ c_{23} (  \nz{3}{3} \wedge \nz{3}{3} + 0 )
\\ 
= &
2 c_{01} \nz{0}{3} \wedge \nz{2}{3} 
+  c_{02} \nz{0}{3} \wedge \nz{3}{3} 
+  3 c_{02} \nz{1}{3} \wedge \nz{2}{3} 
+ ( 3 c_{03}+  c_{12})  \nz{1}{3} \wedge \nz{3}{3}  
+ 2 c_{13} \nz{2}{3} \wedge \nz{3}{3}  
\end{align*} 
By solving a homogeneous linear equations, we get 
$\ds  T = c_{03}
(   \nz{0}{3} \wedge \nz{3}{3} - 3 \nz{1}{3} \wedge \nz{2}{3}) + 
c_{23} \nz{2}{3} \wedge \nz{3}{3} $, and   
$\ds  
 \nz{0}{3} \wedge \nz{3}{3} - 3 \nz{1}{3} \wedge \nz{2}{3}
$ is the highest weight vector of $\kS{0}$ and   
$\ds  
\nz{2}{3} \wedge \nz{3}{3}$ is of $\kS{4}$, and we have 
$\ds  \Lambda^2 \kS{3} = \frakS{0} \oplus \frakS{4}$.   
%\end{exam} 
%\begin{exam}\label{exam::f} 
Thus, 
when $n=1$ and the weight=2 case, we add our decomposition and we see 
\begin{align*}
        \CGF{1}{2}{0}{2} 
&= \kS{4}^{Sp(2,\mR)} = \{ \mathbf{0} \}\\
        \CGF{2}{2}{0}{2} 
&= (\Lambda^2 \kS{3})^{Sp(2,\mR)} =  
(\kS{0} \oplus \frakS{4})^{Sp(2,\mR)} = \frakS{0} = \mR  
\end{align*} 
since $\ds  \kS{q} ^ {Sp(2n,\mR)} = \{\mathbf{0}\}$ for $q >2$.     
\end{exam}

\begin{exam} \label{exam::e} 
$\ds  \Lambda^2 \kS{4}$
is a component of 
$\ds  
\cgf{2}{2n}{1}{4}$
in Example \ref{exam::b}.  When $n=1$, we can decompose it into the irreducible
components as the same way above, and get  
$\ds  \Lambda^2 \kS{4} = \frakS{2} \oplus \frakS{6}$. 

When $n=1$, since 
we have the next decomposition 
$\ds  \Lambda ^3 \kS{5} \cong \frakS{9} \oplus
\kS{5} \oplus \frakS{3}$,   
we have no rule of 
negative 4-step descending nor 2-step descending
for $\ds  \Lambda^p \kS{q}$.  
\end{exam}

We would like to see the effect of tensor product. 

\begin{exam} \label{exam::g} 
When weight=6, we know that  
$\ds  \cgf{2}{2n}{1}{6}$ has 3 components
$(3,7)$, $(4,6)$, $(5^2)$, and 
$\ds  \cgf{3}{2n}{1}{6}$ has 3 components
$(3^2,6)$, $(3,4,5)$, $(4^3)$, and 
$\ds  \cgf{4}{2n}{1}{6}$ has 2 components
$(3^3,5)$, $(3^2,4^2)$, and            
$\ds  \cgf{5}{2n}{1}{6}$ has 1 component
$(3^4,4)$, and 
$\ds  \cgf{6}{2n}{1}{6}$ has 1 component
$(3^6)$.  

Assume $n=1$. Concerning with degree 2 cases, since 
$\ds  \kS{3} \otimes \frakS{7} 
= \kS{4} \oplus \cdots \oplus \frakS{10}$ and 
$\ds  \kS{4} \otimes \frakS{6}
= \kS{2} \oplus \cdots \oplus \frakS{10}$, we see that 
$$\ds  [\kS{3} \otimes \frakS{7}, \frakS{0}] =0 \text{\quad
and\quad}  
\ds  [\kS{4} \otimes \frakS{6}, \frakS{0}] =0$$     
On the other hand, we get $\ds  \Lambda^2 \kS{5} =
\kS{0} \oplus 
\kS{4} \oplus 
\kS{8}$, we have 
$$\ds  [\Lambda^2 \kS{5} , \frakS{0}] =1$$     

degree 3: 
(3.1)-case: $\ds  \Lambda^2( \kS{3}) \otimes \frakS{6}  =
( \kS{0} \oplus  \frakS{4}) \otimes  \frakS{6}  
= \kS{6} + ( \frakS{2} + \cdots + \frakS{10})$, and so  
$$\ds  [\Lambda^2( \kS{3}) \otimes \frakS{6} , \frakS{0}] =
0$$  
\begin{align*} \text{
(3.2)-case:}\quad  
   \kS{3} \otimes \frakS{4} \otimes \frakS{5}= &
( \kS{1} + \frakS{3} + \frakS{5} + \frakS{7}) \otimes \frakS{5} \\
=& ( \kS{4} + \frakS{6} ) + 
(\kS{2} + \frakS{4} + \frakS{6} + \frakS{8})
+ (\kS{0} + \frakS{2} + \frakS{4} + \frakS{6} +\frakS{8} +
\kS{10}) 
\\ & 
+ (\kS{2} + \frakS{4} + \frakS{6} +\frakS{8} + \frakS{10}+ \frakS{12}) 
\end{align*} 
we see that 
$$\ds  [\kS{3} \otimes \frakS{4} \otimes  \frakS{5} ,
\kS{0}] =1$$ Since $\ds  \Lambda^3 \frakS{4} = \frakS{2} +
\kS{6}$,      
$$\ds  [\Lambda^3 \kS{4}, \frakS{0}] = 0$$  
degree 4: 
(4.1)-case: $\ds  (\Lambda^3  \kS{3}) \otimes \frakS{5}  =
 \kS{3}  \otimes  \frakS{5}  
= \kS{2} + \cdots +  \frakS{8}$, here we used the property that 
$\ds  \Lambda^k W = \Lambda ^{\dim W-k} W$ in general, and   
$$\ds  [( \Lambda^3 \kS{3}) \otimes \frakS{5} , \frakS{0}] =
0$$  
(4.2)-case: $\ds  (\Lambda^2  \kS{3}) \otimes  (\Lambda^2
\kS{4}) = (\frakS{0} + \frakS{4}) \otimes ( \frakS{2} + \frakS{6}) 
= ( \kS{2} + \frakS{6}) + ( \frakS{2} + \frakS{4} + \frakS{6}) + 
( \kS{2} + \cdots  + \frakS{10}) $,  
$$\ds  [ (\Lambda^2  \kS{3}) \otimes  (\Lambda^2
\kS{4}), \frakS{0}] = 0 $$  
(5)-case: $\ds  (\Lambda^4  \kS{3}) \otimes \frakS{4}  =
\kS{0} \otimes \frakS{4} = \frakS{4}$ because of $\dim \frakS{3}=4$,
$$\ds  [(\Lambda^4  \kS{3}) \otimes \frakS{4},\frakS{0}]
=0$$  
We add those facts in the table in Example \ref{exam::c}.

\hfil
\begin{tabular}[t]{|c|c|l|c|}
\hline
degree & ref.\# & type &$\dim$ \\\hline
1 & 1 & ($8$)   & 0 \\\hline
2 & 1 & ($3\ 7$) & 0\\
  & 2 & ($4\ 6$) & 0\\
  & 3 & ($5^2 $) & 1\\\hline
\end{tabular}
\hfil 
\begin{tabular}[t]{|c|c|l|c|}
\hline
degree &ref.\# & type &$\dim$ \\\hline 
3 & 1 & ($3^2\ 6$) & 0\\
  & 2 & ($3\ 4\ 5$)& 1\\
  & 3 & ($4^3 $) & 0\\\hline
\end{tabular}
\hfil  
\begin{tabular}[t]{|c|c|l|c|}
\hline
degree &ref.\# & type &$\dim$ \\\hline
4 & 1 & ($3^3\ 5$)& 0\\
  & 2 & ($3^2\ 4^2$)& 0\\\hline
5 & 1 & ($3^4\ 4$)& 0\\\hline
6 & 1 & ($ 3^6$) & 0 \\\hline
\end{tabular}
\hfil 
\end{exam}

\bigskip

We summarize the ideas, so far.  For a given weight $w$ and degree $m$,
we have the complete list of subcomplex of type $\ds 
(k_3,k_3,\ldots, k_s)$ such that $\ds  \sum_{j} (j-2) k_j=w$
and $\ds  \sum_{j} k_j=w$ and 
$$\ds  C^{m}(\kham_{2}^0)_w 
=
\oplus  
\Lambda^{k_3} \kS{3} \otimes 
\Lambda^{k_4} \kS{4} \otimes \cdots 
\Lambda^{k_s} \kS{s} 
$$ 
It is possible to decompose 
$\ds  \Lambda^{p} \kS{q}$ into irreducible subspaces, say, 
$\ds  \Lambda^{p} \kS{q} = \sum \alpha^p_{q r} \frakS{r}$, 
where $\ds  
\alpha^p_{q r} \in \mZ_{\ge 0}$. (But, it is not clear if those  
$\alpha^p_{q r} $ have some ``rule''.)  
Since we have the tensor product formula, called Clebsch-Gordan rule,  
$$ \kS{p} \otimes \frakS{q} =
\kS{|p-q|} \oplus 
\kS{|p-q|+2} \oplus  \cdots 
\kS{p+q}$$ 
we can decompose 
$\ds  
\Lambda^{k_3} \kS{3} \otimes 
\Lambda^{k_4} \kS{4} \otimes \cdots 
\Lambda^{k_s} \kS{s} $ into irreducible components, thus we can
divide 
$\ds  \cgf{m}{2}{1}{w} $ into irreducible components,
and therefore we can pick up the multiplicity of the trivial
representation, and we know the dimension of 
$\ds  \CGF{m}{2}{0}{w}$ as that
multiplicity.  

\begin{thm}
There is a sequence of computer programs which follow the
mathematical story above.  If we input the weight $w$, then we can get 
$\ds  \dim \CGF{m}{2}{0}{w} $ and the
precise contributions of subcomplex, and   
the Euler characteristic number.  If $w$ is big, we may face some trouble
of shortage of memory and so on.  Nevertheless, so for we have a table of 
dimensions of relative cochain complex until $w \le 20$.  
The table below means the horizontal direction is the degree of cochain
complex and the descending vertical direction means increasing weight.

\medskip

\input{thm1.tex}

\medskip

Here, $\chi$ means the alternating sum of the dimension of relative
cochain complex, including 0-dimensional cochain complex $\mR$. Thus, 
$\chi \ne 1$ means there is non-trivial Betti number.  
\end{thm}
\bigskip

%%%%%%%%%%%%%%%%%%%%%%%%%%%%%%%%%%%%%%%%%%%%%%%%%%%%%%%%

\subsection{How to use computer in order to get 
the dimension of relative cochain complex (brief summary)} 
\label{howto:one}

\begin{enumerate}
\item  For a given weight $w$, edit and fix the weight $w$ in 
\texttt{../gkf\_act-1.mpl} and run this  maple script in order to get possible
direct summands of 
$\ds  \cgf{dp}{2}{1}{w}$ for each degree $dp$.  

%\texttt{\% cd //METOKI/GF\_$w$}

\texttt{\% maple ../gkf\_act-1.mpl > OUT-$w$}

\item 
Using the output file 
\texttt{OUT-$w$}, we run the next perl script 

\texttt{\% perl ../gkf\_act-2.prl OUT-$w$}

After this job, we have a couple of output 
\texttt{range-$w$} and files \texttt{cases\_w\_dp.txt}.  

Also, we run 

\texttt{\% perl ../mk-tex-gkf.prl OUT-$w$}

we get a \TeX{}-file  \texttt{gkf\_$w$.tex}.   
Compiling this, we get a table of degree, reference numbers and decompose
type.   

\item 
Key issue here is: First run perl, 

\texttt{\% perl ../Unify/eating.prl gkf\_w.tex}

Then we get a output file \texttt{esa\_w.txt}

In \texttt{../Unify/sophi.mpl}, (1) edit and fix the weight $w$, and (2)
fix path to \texttt{range-w} and \texttt{esa\_w.txt}.  Then run   

\texttt{\% maple -q  ../Unify/sophi.mpl}
    
Then we see contribution of each type and the multiplicity of each cochain
complex.  
\end{enumerate}

%%%%%%%%%%%%%%%%%%%%%%%%%%%%%%%%%%%%%%%%%%%%%%%%%%%%%%%%

\section{Getting concrete bases and matrix representation of $\myd$}
Even though we know all the dimensions of subcomplex of relative cochain complex 
$\ds  \cgf{m}{2n}{1}{w}$, it is not
enough to investigate the coboundary operator $\myd$ itself.  We have to know
concrete bases of 
$\ds  \CGF{m}{2n}{0}{w} $ and 
$\ds  \CGF{m+1}{2n}{0}{w}$, and the matrix
representation of $\myd$ with respect to those bases. 
In this section, we only deal with $n=1$ and $\ds  M=\mR^2$.

We follow the requirements for relative
cochain $\sigma$. The first one is $\ds  i_{\xi} \sigma =0$ 
($\forall \xi \in \fraksp(2,\mR)$). This is very easy to check $\sigma$
contains $\ds  \kS{2}$ and indeed, we    
already have omitted $\ds  \kS{2}$.    
The other thing is  
$\ds  i_{\xi} \myd \sigma =0$  
($\forall \xi \in \fraksp(2,\mR)$).  
Since we have already observed that $\ds  i_{\xi}\myd = -
\xi\cdot $, we get 
\begin{align*}
\left(i_{H}\circ \myd\right) \nz{r}{R} &= - H \cdot \nz{r}{R} =  
 - (2r- R) \nz{r}{R} \\ 
\left(i_{X}\circ \myd\right) \nz{r}{R} &=  
- X \cdot \nz{r}{R} = 
 - (R-r) \nz{r+1}{R}  
\\
\left(i_{Y}\circ \myd\right) \nz{r}{R} &= - Y \cdot \nz{r}{R} 
- r \nz{r-1}{R} 
\end{align*} 
for each generators of 1-cochain complex. 
$\myd$ is a skew-derivation of degree $+1$ and  
$$
\myd \nz{r}{R} = - \frac{r!(R-r)!}{2}  
\sum_{ a+b = 1+r, A+B=2+R }
\begin{vmatrix} a & b\\
A & B\end{vmatrix}
\dfrac{ \nz{a}{A} }{a!  (A-a)! } 
\wedge 
\dfrac{ \nz{b}{B} }{b!   (B-b)!} 
$$
where 
$ 0 \le r \le R$, $0\le a\le A$ , $0\le b \le B$, $A\ge 2$, $B\ge 2$.

\subsection{An easy exmaple} 
\label{exam::w::six}
Here,  we show a small practice of getting concrete basis
of relative cochain complex of the case of $n=1$ and the weight 6.  As
shown in Example \ref{exam::g} by some help of $Sp(2,\mR)$-theory, 
we have already known that 
$\ds  \dim \CGFF{2}{6} = \dim \CGFF{3}{6} =1 $ 
and the others are 0-dimensional, and furthermore 
the corresponding basis lives in 
$\ds  \Lambda^2 \kS{5}$ when 
$\ds  \CGFF{2}{6} $, and 
lives in 
$\ds  \kS{3} \wedge \frakS{4} \wedge \frakS{5} $ when 
$\ds \CGFF{3}{6}$.  
Let $\ds  \sigma = \sum_{0\le i< j \le 5} c_{i,j} \nz{i}{5}
\wedge \nz{j}{5}$, and solve the equations $\ds  i_{\xi} \circ
\myd \sigma = 0$ for $\xi \in \{ X, Y, H\}$, 
which are a basis of $\fraksp(2,\mR)$.  
Then we have 
% var_num:=15;eqn_num:=37;term_num:=3;
$$
\sigma = c_{2,3} \left(
10 \nz{0}{5}\wedge \nz{5}{5}
- 2\nz{1}{5}\wedge \nz{4}{5}
+
\nz{2}{5}\wedge \nz{3}{5}\right)$$
and putting $\ds  c_{2,3}=1$, we have $\ds 
\sigma_1$, a basis of 
$\ds  \CGFF{2}{6}$.  

For 3-cochain, 
let $\ds  \tau  = 
\sum_{ i=0}^{3} 
\sum_{ j=0}^{4} 
\sum_{ k=0}^{5} 
c_{i,j,k} 
\nz{i}{3} \wedge \nz{j}{4}\wedge \nz{k}{5}
$, and again solve the equations $\ds  i_{\xi} \circ
\myd \sigma = 0$ for $\xi \in \{ X, Y, H\}$.  
Then we have 
% var_num:=120;eqn_num:=337;term_num:=18;
\begin{align*}
\tau  & = c_{2,3,1} \left( \nz{2}{3}\wedge \nz{3}{4}\wedge \nz{1}{5}
+8 \nz{1}{3}\wedge \nz{4}{4}\wedge \nz{1}{5} 
-5/2 \nz{1}{3}\wedge \nz{3}{4}\wedge \nz{2}{5}
+\nz{1}{3}\wedge \nz{2}{4}\wedge \nz{3}{5}
-6 \nz{0}{3}\wedge \nz{2}{4}\wedge \nz{4}{5}\right. \\& 
+15 \nz{0}{3}\wedge \nz{1}{4}\wedge \nz{5}{5} 
+\nz{1}{3}\wedge \nz{1}{4}\wedge \nz{4}{5}
-20 \nz{1}{3}\wedge \nz{0}{4}\wedge \nz{5}{5} 
-6 \nz{0}{3}\wedge \nz{4}{4}\wedge \nz{2}{5}
+9/2 \nz{0}{3}\wedge \nz{3}{4}\wedge \nz{3}{5}\\& 
+8 \nz{2}{3}\wedge \nz{0}{4}\wedge \nz{4}{5}
-5/2 \nz{2}{3}\wedge \nz{1}{4}\wedge \nz{3}{5} 
+\nz{2}{3}\wedge \nz{2}{4}\wedge \nz{2}{5}
-20 \nz{2}{3}\wedge \nz{4}{4}\wedge \nz{0}{5}
-6 \nz{3}{3}\wedge \nz{0}{4}\wedge \nz{3}{5}\\& \left.
+9/2 \nz{3}{3}\wedge \nz{1}{4}\wedge \nz{2}{5} 
-6 \nz{3}{3}\wedge \nz{2}{4}\wedge \nz{1}{5}
+15 \nz{3}{3}\wedge \nz{3}{4}\wedge \nz{0}{5}\right) 
\end{align*}
and putting $\ds  c_{2,3,1}=4$, we have a basic vector 
$\tau_1$ in 
$\ds\CGFF{3}{6}$.  We see that  
$$ \myd \sigma_1 = \tau_1$$
thus, we get 
$\ds  \HGFF{2}{6}  = \{ \mathbf{0}\}$, and so we have 
$\ds  \HGFF{3}{6}  = \{ \mathbf{0}\}$. Therefore, the Euler
characteristic number of the alternating sum of the Betti numbers, which
we 
include the trivial 0-th Betti number 1,  is 1 as we have seen before.

\subsection{More complicated example \label{subsec:exam}}
By taking the example of weight = 8, 
we emphasize that how getting concrete basis is tough job.  
Later, we stress that 
our job sequences, which complete those jobs automatically, are how
useful in this discussion.

%\tablelasttail{\hline}
\hfil \begin{tabular}[t]{cr|c r }%\hline
\hline deg  & ref.\# & type & $\dim$ \\
\hline 1 &  1 & ($10$)\\
\hline 2 &  1 & ($3\ 9$)\\
2 &  2 & ($4\ 8$)\\
2 &  3 & ($5\ 7$)\\
2 &  4 & ($6^2$)\\
\hline 3 &  1 & ($3^2\ 8$)\\
3 &  2 & ($3\ 4\ 7$) & 1\\
3 &  3 & ($3\ 5\ 6$) & 1\\
3 &  4 & ($4^2\ 6$) & 1\\
3 &  5 & ($4\ 5^2$) & 1\\
\hline
\end{tabular}
\hfil 
% 
%\tablefirsthead{\multicolumn{4}{c}{weight= 8}\\\hline}
%\tablehead{\multicolumn{4}{c}{weight= 8}\\\hline}
%\tablelasttail{\hline}
\begin{tabular}[t]{cr|c r }%\hline
\hline deg  & ref.\# & type & $\dim$ \\ 
\hline 4 &  1 & ($3^3\ 7$)\\
4 &  2 & ($3^2\ 4\ 6$) & 1\\
4 &  3 & ($3^2\ 5^2$) & 2\\
4 &  4 & ($3\ 4^2\ 5$) & 2\\
4 &  5 & ($4^4$)\\
\hline 5 &  1 & ($3^4\ 6$)\\
5 &  2 & ($3^3\ 4\ 5$) & 1\\
5 &  3 & ($3^2\ 4^3$)\\
\hline 6 &  1 & ($3^4\ 4^2$)\\
\hline
\end{tabular}
\hfil 

The cochain complex of deg = 3 is 4-dim, that of deg = 4 is
5-dim, and that of deg = 5
is 1-dimensional. 

Here, we will see the whole process of getting a basis of each cochain
complex.  

\begin{enumerate}
\item  3-cochains 

\begin{enumerate}
\item  type ($3, 4, 7$)-case:  
The candidate is $\ds  
\sigma_1 = 
\sum_{i=0}^ {3} 
\sum_{j=0}^ {4} 
\sum_{k=0}^ {7} c_{ijk} \nz{i}{3} \wedge \nz{j}{4} \wedge \nz{k}{7}$
with 160 unknown variables $c_{ijk}$. 
From the three conditions $\ds  i_{\xi} \circ \myd \sigma=0$
for $\xi = X, Y, H$, we have a homogeneous linear equations ( the number
of equations is 455). Solving these equations, we have one undetermined
parameter $c_{007}$, 
and we put it to be 1,   
and get a 1-dimensional basis:  

%var_num:=160;eqn_num:=455;term_num:=20;

\item  type ($3, 5, 6$)-case: The method is completely same and get a
basis $\ds  \sigma_2$ by 
putting $c_{016} = 1$.  
% var_num:=168;eqn_num:=477;term_num:=22;

\item  type ($4^2, 6$)-case: Here we have to be careful for dealing with 
exterior product elements as below. 
The candidate is $\ds  
\sigma = 
\sum_{0\le i< j\le 4} 
\sum_{k=0}^{6} c_{ijk} \nz{i}{4} \wedge \nz{j}{4} \wedge \nz{k}{6}$
with 70 unknown variables $c_{ijk}$. Putting $\ds  c_{016}=1$,
we have a basis $\ds  \sigma_3$. 
% var_num:=70;eqn_num:=199;term_num:=10;

\item  type ($4, 5^2$)-case: Here again we have to be careful for dealing with 
exterior product elements as before. 
The candidate is $\ds  
\sigma_4 = 
\sum_{i=0}^{4} 
\sum_{0\le j< k\le 5} 
 c_{ijk} \nz{i}{4}\wedge \nz{j}{5} \wedge \nz{k}{5}$
with 75 unknown variables $c_{ijk}$. Putting $c_{025}=1$, we have a
basis
% var_num:=75;eqn_num:=213;term_num:=11; 
\end{enumerate}

\input in_dp3_w8.tex

Those $\ds  \sigma_1, \sigma_2, \sigma_3, \sigma_4$ are a basis
of 
$\ds  \CGFF{3}{8}$.

\item  4-cochains 
\begin{enumerate}
\item  type ($3^2, 4, 6$)-case: 
$ \ds  c_{0206} \tau_1 $, 
where $\tau_1$ is below:
%var_num:=210;eqn_num:=599;term_num:=24; 

\item  type ($3^2, 5^2$)-case: The general solution is 
$ \ds  c_{0314} \tau_2 + c_{0314} \tau_3$, 
where $\tau_2, \tau_3$ are below:
% var_num:=90;eqn_num:=253;term_num:=table([(1)=13}{(2)=13};

\item  type ($3, 4^2, 5$)-case: 
The general solution is $ \ds 
c_{0125} \tau_4 + c_{0035} \tau_5$, where $\tau_4, \tau_5$ are below:
% 
%%%%%%%%%%%%%%%%%%%%%%%%%%%%%%%%%%%%%%%%%%%%%%%%%%%%%%%%%%%
% var_num:=240;eqn_num:=683;term_num:=table([(1)=32}{(2)=32};
%
\end{enumerate}

\input in_dp4_w8.tex

\item  5-cochain 

type ($3^3\ 4\ 5$)-case: 120 unknown variables and 
solving the linear 339-equations, and putting $c_{01314}=1$, we have a
basis:
% var_num:=120;eqn_num:=339;term_num:=18;

\input in_dp5_w8.tex

\end{enumerate}

With respect to those concrete bases of relative cochain
complex, the coboundary operator $\myd$ between degree 4 and 5,    
is of the form $\ds  \myd \tau_j = b_j \rho \quad ( j=1..5)$
and comparing both sides, we get a matrix expression 
$\ds  B := \begin{pmatrix} b_1 & b_2 & b_3 & b_4 & b_5 \end{pmatrix}
= \begin{pmatrix} -4 & -3 & -3  & -9  & 6 \end{pmatrix}$. 

About the coboundary operator $\myd$ between degree 3 and 4, 
$\myd \sigma_j = \sum_{i=1}^{5} a_j ^ {i} \tau_i$ ($j=1,\ldots, 4$), For
each $i$, comparing both sides again, we get a 
matrix expression 
$$\ds  A = \left( a^{j}_{i} \right)  
= \begin{pmatrix} -6  & -12 & \frac{9}{4}  & 0   \\
0 & -27 & 0 & 1 \\
0 & 59 & 0 & 3 \\
5 & 4 &  \frac{3}{4} & 22 \\
\frac{7}{2} & 14 & \frac{21}{8}  & 35 
\end{pmatrix}$$ 
As expected, we see that $\ds  B A = O$, which corresponds to $\myd \circ
\myd = 0$. Since $\rank A=4$, $\rank B=1$, thus we have
$\ds  \HGFF{3}{8} = \{ \mathbf{0} \}$, 
$\ds  \HGFF{4}{8} = \{ \mathbf{0} \}$, 
$\ds  \HGFF{5}{8} = \{ \mathbf{0} \}$.  These mean that except  
$\ds  \HGFF{0}{8} = \mR $,  the other cohomologies are zero and 
the Euler characteristic number is 1.  

\subsection{How tuning computer to subscribe our job}
We here show our typical strategy for getting a concrete basis of
relative cochain complex. Take an example with weight 8, which we have
seen above.  We follow the above discussion in the case of type
$(3,5,6)$ and type $\ds  (4^2,6)$ in order to emphasize some
difference.  
We need to distinguish scalar and vectors in general. For
that purpose, we use \texttt{difforms} package.  
In the type ($3,5,7$) of degree 3 with weight 8. a  cochain in general is  
$$ \sigma = \sum_{i_1=0}^ {3} \sum_{i_2=0}^ {5} \sum_{i_3=0}^ {7} 
c_{i_1, i_2, i_3} \nz{i_1}{3}  \otimes
\nz{i_2}{5} \otimes \nz{i_3}{7} $$(in this case, $\wedge$ and $\otimes$
have the same meaning), and we declare $c_{i_1, i_2, i_3}$ are scalar. 
For that purpose
we prepare \texttt{in\_8\_3-3\_a.txt}. 
\begin{quote}\small
\begin{verbatim}
uke31 := NULL: uke32 := NULL: 
for i1 from 0 to d3 do
for i2 from 0 to d5 do
for i3 from 0 to d6 do
   uke31 := uke31, cat(c3, "_",i1, "_",i2, "_",i3) = const ;
   uke32 := uke32, cat(c3, "_",i1, "_",i2, "_",i3) ;
od od od:
defform(uke31):
\end{verbatim}
\end{quote}
To construct $\sigma$ above, we prepare \texttt{in\_8\_3-3\_b.txt}.  
\begin{quote}\small
\begin{verbatim}
# defform( A3 = 3): A3 :=0:
for i1 from 0 to d3 do
for i2 from 0 to d5 do
for i3 from 0 to d6 do
A3 := A3 + cat(c3, "_",i1, "_",i2, "_",i3) * 
&^ ( z[i1,d3-i1], z[i2,d5-i2], z[i3,d6-i3] ): 
od od od:
\end{verbatim}
\end{quote}

In the type 
($4^2,6$) of degree 3 with weight 8.  a general cochain is  
$$ \sigma = \sum_{i_1 =0}^{4}\sum_{i_2 =0}^{4}
\sum_{i_3=0}^{6} c_{i_1,i_2,i_3} 
\nz{i_1}{4} \wedge 
\nz{i_2}{4}  \wedge  \nz{i_3}{6}, 
\quad\text{where}\   c_{i_1,i_2,i_3} \ \text{are skewsymmetric in 3
indices.  }
$$
To avoid complicated requirement, we restrict range of indices as
$$ \sigma = \sum_{0 \le i_1 < i_2 \le 4}
\sum_{i_3=0}^{6} c_{i_1,i_2,i_3} 
\nz{i_1}{4} \wedge 
\nz{i_2}{4}  \otimes   \nz{i_3}{6}, 
$$ (we need some scale factor, in the case above, 2-times, but we ignore
them here after.)

This time, 
a declaration that coefficients are constant, is  
\texttt{difforms}, 
we prepare \texttt{in\_8\_3-4\_a.txt}. 
\begin{quote}\small
\begin{verbatim}
uke41 := NULL: uke42 := NULL: 
for i1 from 0 to d4 do
for i2 from 1+i1 to d4 do
for i3 from 0 to d6 do
   uke41 := uke41, cat(c4, "_",i1, "_",i2, "_",i3) = const ;
   uke42 := uke42, cat(c4, "_",i1, "_",i2, "_",i3) ;
od od od:
defform(uke41):
\end{verbatim}
\end{quote}
To construct $\sigma$ above, we prepare \texttt{in\_8\_3-4\_b.txt}.  
\begin{quote}\small
\begin{verbatim}
# defform( A4 = 3): A4 :=0:
for i1 from 0 to d4 do
for i2 from 1+i1 to d4 do
for i3 from 0 to d6 do
A4 := A4 + cat(c4, "_",i1, "_",i2, "_",i3) * 
&^ ( z[i1,d4-i1], z[i2,d4-i2], z[i3,d6-i3] ): 
od od od:
\end{verbatim}
\end{quote}
%\end{me}

To pick up ``coefficients'' of $\ds  i_{\xi} \myd(\sigma) $ ($\forall \xi\in
\fraksp(2,\mR)$) with respect ``some generators'', 
we prepare \texttt{in\_8\_3-4\_c.txt}.  
\begin{quote}\small
\begin{verbatim}
for i1 from 0 to d4 do
for i2 from 1+i1 to d4 do
for i3 from 0 to d6 do
ukez := ukez, seq(
mInnd([z[i1,d4-i1],z[i2,d4-i2],z[i3,d6-i3]],
                           W_W[s]),s=0..d0);
od od od:
\end{verbatim}
\end{quote}
%\end{me}

For the main purpose to solve  $\ds  i_{\xi} \myd(\sigma)=0 $ ($\forall \xi\in
\fraksp(2,\mR)$) and to determine $\sigma$ and $\myd(\sigma)$, we have a core
Maple script
\texttt{action\_new.mpl}, 
we omit it because it will occupy one and a half pages.    

\begin{myRemark}

When  $ w \ge 16$, we encounter some
trouble of kernel panic in processing the above steps.  Main problem seems
shortage of CPU memory. To recover
this trouble, we change our process slightly.   
A big guarantee is {\em linearity} of $\myd$. Namely,  
suppose for a given  ``huge'' cochain $A$  we cannot compute
$\myd(A)$. Then, we divide $A$ into ``small'' pieces like 
$\ds  A = a_1 + a_2 + \cdots$.   
Instead of handling $\myd(A)$ itself, we manipulate  
$\ds  \myd
( a_i)$, and by $\ds  \myd (a_1) + \myd (a_2) + \cdots$, we get the
whole 
$\ds  \myd (A)$. 
\end{myRemark}

\subsection{Sequence of computer process to get concrete basis of
cochain complex (brief summary)}  

\begin{itemize}
\item[1.] The same (1) of \S \ref{howto:one}. 
\item[2.] The same (2) of \S \ref{howto:one}. 

\item[3.]
We prepare a plain text file 
\texttt{in\_dp-new.txt}, 
this is a prototype of cochain of
degree $dp$, independent of $w$.  
In which we prepare coefficients of degree $dp$
cochain, and ready to handle to 
defform() of difforms Maple package, and 
(Of course, the range of $dp$ depends on $w$.)

We manipulate this \texttt{in\_dp-new.txt} and 
\texttt{cases\_w\_dp.txt} in step 2 by the perl script
\texttt{gkf\_act-3.prl}.

\texttt{\% perl ../gkf\_act-3.prl}

Out put file are 
\texttt{out\_w\_dp-ref\#}.

\item[4.] \textbf{Important!} 
For \texttt{out\_w\_dp-ref\#}, we have to revise it in order to guarantee the
skew-symmetry. For instance, the right side is desired form: 
\begin{quote}
\begin{tabular}{l|l}
\texttt{for i1 from 0 to d4 do } & \texttt{for i1 from 0 to d4 do } \\
\texttt{for i2 from 0 to d4 do } & \texttt{for i2 from \underline{1+i1} to d4 do} \end{tabular} 
\end{quote} 

\medskip

For that purpose, we prepare a 
perl script  \texttt{gkf\_act-4.prl} and run  for all files: 

\texttt{\% ../gkf\_act-4.prl out\_w\_dp-ref\#} 

Output are 
\texttt{out\_w\_dp-ref\#.txt}.   

\item[5.] For the revised 
\texttt{out\_w\_dp-ref\#.txt}, we apply \texttt{gkf\_act-5.prl}.  

\texttt{\% ../gkf\_act-5.prl out\_w\_dp-ref\#.txt}

Output are 
\texttt{in\_w\_dp-ref\#\_[abc].txt}.   In \S \ref{subsec:exam}, we
encountered with \texttt{in\_8\_3-4\_[abc].txt}.   

\item[6.]
We prepare a maple script 
\texttt{action\_new.mpl}, in which for fixed weight $w$ and  
degree  $dp$, get a basis of cocycles and 
$d$-image of them for each direct sum component labeled by ref\#.  

We omitted to explain about this key job in \S \ref{subsec:exam}.

\end{itemize}

\section{Main result} 
While weight is less that 12, 
the whole cohomology groups were studied in \cite{KOT:MORITA}.

\begin{thm} 
Using those simple but powerful tricks explained above, we could finish all
computations for 
weight from 12 to 18 of $\ds \HGF{\bullet}{2}{0}{wt}$.  We will show the results a table below. The    
abbreviations in the table mean 
degree $k$ is  $\texttt{Sp}(2,\mR)$-invariant 
cochain complex $\ds  \text{C}^k$,  $\dim$ is $\dim \text{C}^k$, and   
$\rank(d)$ is the rank of $ \ds  d : \text{C}^k \longrightarrow
\text{C}^{1+k}$.  
% $\ds  \CGFF{\bullet}{12}$
\begin{center}
\begin{tabular}[t]{|c|c|*{10}{r}|}
\hline 
weight & degree  & 1& 2& 3& 4 & 5& 6 &7& 8 &9 & 10 \\\hline
12 & $\dim$  & 0& 0& 8&23 &22&13 &5& 0 & 0 & 0 \\ 
   & $\rank(d)$ & 0& 0& 8&14 & 8& 5 &0& 0 & 0 & 0\\ 
   & $Betti \#$ &  0& 0& 0& 1 & 0& 0 &0& 0& 0 & 0 \\
\hline 
14 & $\dim$  & 0& 1& 6&31 &71&58 &15& 2 &1 & 0\\ 
   & $\rank(d)$ & 0& 1& 5&26 & 44& 14 &1& 1 & 0 & 0\\ 
   & $Betti \#$ & 0& 0& 0& 0 & 1& 0 &0& 0&0 & 0 \\
\hline 
16 & $\dim$  & 0& 0& 12 &61 &126 & 147 &95& 24 &0 & 0\\ 
   & $\rank(d)$ & 0& 0& 12&49 & 77& 70 &24& 0 & 0 & 0\\ 
   & $Betti \#$ & 0& 0& 0& 0 & 0& 0 &1& 0&0 & 0 \\
\hline 
18 & $\dim$  & 0& 1& 10 &80 &262 & 380 &268& 100 &21 & 1\\ 
   & $\rank(d)$ & 0& 1& 9&71 & 191& 188 &80& 20 & 1 & 0\\ 
   & $Betti \#$ & 0& 0& 0& 0 & 0& 1 &0& 0&0 & 0 \\
\hline 
\end{tabular} 
\end{center}

\bigskip

\end{thm}

\begin{myRemark}
The Euler characteristic number is the alternating sum of $dim$ of cochain complexes or Betti
numbers, including 0-dimensional. 
The tables above show that 
$\chi(\text{weight}=12) = 2$, 
$\chi(\text{weight}=14) = 0$, 
$\chi(\text{weight}=16) = 0$, and  
$\chi(\text{weight}=18) = 2$.  

When the weight = 20, 
we have the complete list of degree of all relative cochain complexes as
below:

\medskip

\hfil 
\begin{tabular} { |c|*{11}{r}|}
\hline
degree & 1 & 2 & 3 & 4 & 5 &  6 & 7 & 8 & 9 & 10 & 11 \\\hline 
$\dim$ & 
0  & 0  & 17 & 124& 423& 791& 801 & 414 & 96 & 9 & 1 \\\hline
\end{tabular}

\bigskip

and the Euler characteristic number, i.e.,  
the alternating sum including 0-dimensional, is 1.  It will be interesting to check if there are
non-trivial cohomology groups with weight 20.  
\end{myRemark}

\begin{myRemark}
        The above remark was noted on October, 2012. Recently (February 2014)
        we finished computing the Betti numbers when weight 20, using
        Gr\"obner Basis theory for linear polynomials.  Our result  is that 
        the Betti numbers for weight 20 
        are all trivial, except 0-dimensional.   
\end{myRemark}

%\cite{metoki:shinya}

%\nocite{morita:text}
\nocite{okada:text}
%\nocite{M:Takamura}
\nocite{fulton:harris}
\nocite{goodman:wallach}

\bibliographystyle{plain}
\bibliography{km_refs}

\end{document}

%% file: kmme11.tex
\long\def\kmcomment#1{}
\usepackage{amsmath,amssymb}

\usepackage[a4paper,truedimen,scale={.80,.85},centering]{geometry} 

{\list{}{\leftmargin=#1 \topsep=0pt \parsep=0pt \itemsep=0pt}\item[]}%
{\endlist}
\RequirePackage{verbatim}
{\noindent\quote\endgraf\verbatim}%
{\endverbatim\endquote\endgraf\medskip}%
\makeatletter
\newenvironment{qverb*}%
{\noindent\quote\endgraf\@nameuse{verbatim*}}%
{\@nameuse{endverbatim*}\endquote\endgraf\medskip}%
\makeatother

%% file: thm1.tex
%\input kmmj11
%\begin{document}
\hfil 
\begin{tabular}{|r|*{14}{r}|c|}
\hline
$w$ & 1& 2& 3& 4& 5& 6& 7& 8& 9& 10& 11& 12 & 13 & 14 & $\chi$ \\\hline
 2  & 0& 1&  &  &  &  &  &  &  &   &   &&&& $2$ \\\hline 
 4  & 0& 0& 1& 1&  &  &  &  &  &   &   &&&& $1$ \\\hline 
 6  & 0& 1& 1& 0& 0&  &  &  &  &   &   &&&& $1$ \\\hline 
 8  & 0& 0& 4& 5& 1& 0&  &  &  &   &   &&&& $1$ \\\hline 
10  & 0& 1& 3& 9&12& 4& 0&  &  &   &   &&&& $0$ \\\hline 
12  & 0& 0& 8&23&22&13& 5& 0&  &   &   &&&& $2$ \\\hline 
14  & 0& 1& 6&31&71&58&15& 2& 1&   &   &&&& $0$ \\\hline 
16  & 0& 0&12&61&126&147&95&24& 0&   &   &&&& $0$ \\\hline 
18  & 0& 1&10&80&262&380&268&100&21&1&   &&&& $2$ \\\hline 
20  & 0& 0&17&124&423&791&801&414&96&9&1 &&&& $1$ \\\hline 
22  &0&1&14&163&738&1586&1874&1276&479&82&3&&&&$1$  \\\hline 
24  &0&0&23&229&1091&2897&4281&3534&1628&408&49&1&&&$-2$\\\hline 
26  &0&1&19&285&1722&5102&8613&8735&5222&1703&266&19&1&&$3$ \\\hline 
28 &0&0&29&385&2428&8465&16905&19930&14133&5981&1408&144&2& &$1$ \\\hline 
30 &0&1&25&468&3541&13661&30687&42291&35986&18457&5431&855&63&0&$1$\\\hline 
% x  & 1& 2& 3& 4& 5& 6  & 7& 8& 9& 10& 11& $\chi$ \\\hline 
%    &  &  &  &  &  &  &  &  &  &   &   & $\chi$ \\\hline 
\end{tabular}

\medskip
%\end{document} 

%% file: in_dp3_w8.tex
\begin{align*}
\sigma_1=& \nz{0}{3}\wedge \nz{0}{4}\wedge \nz{7}{7}
-4\nz{0}{3}\wedge \nz{1}{4}\wedge \nz{6}{7}
+6\nz{0}{3}\wedge \nz{2}{4}\wedge \nz{5}{7}
-4\nz{0}{3}\wedge \nz{3}{4}\wedge \nz{4}{7}
+\nz{0}{3}\wedge \nz{4}{4}\wedge \nz{3}{7}
\\&
-3\nz{1}{3}\wedge \nz{0}{4}\wedge \nz{6}{7}
+12\nz{1}{3}\wedge \nz{1}{4}\wedge \nz{5}{7}
-18\nz{1}{3}\wedge \nz{2}{4}\wedge \nz{4}{7}
+12\nz{1}{3}\wedge \nz{3}{4}\wedge \nz{3}{7}
-3\nz{1}{3}\wedge \nz{4}{4}\wedge \nz{2}{7}
\\&
+3\nz{2}{3}\wedge \nz{0}{4}\wedge \nz{5}{7}
-12\nz{2}{3}\wedge \nz{1}{4}\wedge \nz{4}{7}
+18\nz{2}{3}\wedge \nz{2}{4}\wedge \nz{3}{7}
-12\nz{2}{3}\wedge \nz{3}{4}\wedge \nz{2}{7}
+3\nz{2}{3}\wedge \nz{4}{4}\wedge \nz{1}{7}
\\&
-\nz{3}{3}\wedge \nz{0}{4}\wedge \nz{4}{7}
+4\nz{3}{3}\wedge \nz{1}{4}\wedge \nz{3}{7}
-6\nz{3}{3}\wedge \nz{2}{4}\wedge \nz{2}{7}
+4\nz{3}{3}\wedge \nz{3}{4}\wedge \nz{1}{7}
-\nz{3}{3}\wedge \nz{4}{4}\wedge \nz{0}{7}
\\
\sigma_2=&\nz{0}{3}\wedge \nz{1}{5}\wedge \nz{6}{6}
-4\nz{0}{3}\wedge \nz{2}{5}\wedge \nz{5}{6}
+6\nz{0}{3}\wedge \nz{3}{5}\wedge \nz{4}{6}
-4\nz{0}{3}\wedge \nz{4}{5}\wedge \nz{3}{6}
+\nz{0}{3}\wedge \nz{5}{5}\wedge \nz{2}{6}
\\&
-\nz{1}{3}\wedge \nz{0}{5}\wedge \nz{6}{6}
+2\nz{1}{3}\wedge \nz{1}{5}\wedge \nz{5}{6}
+2\nz{1}{3}\wedge \nz{2}{5}\wedge \nz{4}{6}
-8\nz{1}{3}\wedge \nz{3}{5}\wedge \nz{3}{6}
+7\nz{1}{3}\wedge \nz{4}{5}\wedge \nz{2}{6}
\\&
-2\nz{1}{3}\wedge \nz{5}{5}\wedge \nz{1}{6}
+2\nz{2}{3}\wedge \nz{0}{5}\wedge \nz{5}{6}
-7\nz{2}{3}\wedge \nz{1}{5}\wedge \nz{4}{6}
+8\nz{2}{3}\wedge \nz{2}{5}\wedge \nz{3}{6}
-2\nz{2}{3}\wedge \nz{3}{5}\wedge \nz{2}{6}
\\&
-2\nz{2}{3}\wedge \nz{4}{5}\wedge \nz{1}{6}
+\nz{2}{3}\wedge \nz{5}{5}\wedge \nz{0}{6}
-\nz{3}{3}\wedge \nz{0}{5}\wedge \nz{4}{6}
+4\nz{3}{3}\wedge \nz{1}{5}\wedge \nz{3}{6}
-6\nz{3}{3}\wedge \nz{2}{5}\wedge \nz{2}{6}
\\&
+4\nz{3}{3}\wedge \nz{3}{5}\wedge \nz{1}{6}
-\nz{3}{3}\wedge \nz{4}{5}\wedge \nz{0}{6}
\\
\sigma_3=&\nz{0}{4}\wedge \nz{1}{4}\wedge \nz{6}{6}
-3\nz{0}{4}\wedge \nz{2}{4}\wedge \nz{5}{6}
+3\nz{0}{4}\wedge \nz{3}{4}\wedge \nz{4}{6}
-\nz{0}{4}\wedge \nz{4}{4}\wedge \nz{3}{6}
+6\nz{1}{4}\wedge \nz{2}{4}\wedge \nz{4}{6}
\\&
-8\nz{1}{4}\wedge \nz{3}{4}\wedge \nz{3}{6}
+3\nz{1}{4}\wedge \nz{4}{4}\wedge \nz{2}{6}
+6\nz{2}{4}\wedge \nz{3}{4}\wedge \nz{2}{6}
-3\nz{2}{4}\wedge \nz{4}{4}\wedge \nz{1}{6}
+\nz{3}{4}\wedge \nz{4}{4}\wedge \nz{0}{6}
\\ 
\sigma_4=& \nz{0}{4}\wedge \nz{2}{5}\wedge \nz{5}{5}
-3\nz{0}{4}\wedge \nz{3}{5}\wedge \nz{4}{5}
-2\nz{1}{4}\wedge \nz{1}{5}\wedge \nz{5}{5}
+4\nz{1}{4}\wedge \nz{2}{5}\wedge \nz{4}{5}
+\nz{2}{4}\wedge \nz{0}{5}\wedge \nz{5}{5}
\\&
+\nz{2}{4}\wedge \nz{1}{5}\wedge \nz{4}{5}
-8\nz{2}{4}\wedge \nz{2}{5}\wedge \nz{3}{5}
-2\nz{3}{4}\wedge \nz{0}{5}\wedge \nz{4}{5}
+4\nz{3}{4}\wedge \nz{1}{5}\wedge \nz{3}{5}
+\nz{4}{4}\wedge \nz{0}{5}\wedge \nz{3}{5}
\\&
-3\nz{4}{4}\wedge \nz{1}{5}\wedge \nz{2}{5}
\end{align*}

%% file: in_dp4_w8.tex
\begin{align*}
\tau_1=& -2\nz{0}{3}\wedge \nz{1}{3}\wedge \nz{1}{4}\wedge \nz{6}{6}
+6\nz{0}{3}\wedge \nz{1}{3}\wedge \nz{2}{4}\wedge \nz{5}{6}
-6\nz{0}{3}\wedge \nz{1}{3}\wedge \nz{3}{4}\wedge \nz{4}{6}
+2\nz{0}{3}\wedge \nz{1}{3}\wedge \nz{4}{4}\wedge \nz{3}{6}
\\&
+\nz{0}{3}\wedge \nz{2}{3}\wedge \nz{0}{4}\wedge \nz{6}{6}
-6\nz{0}{3}\wedge \nz{2}{3}\wedge \nz{2}{4}\wedge \nz{4}{6}
+8\nz{0}{3}\wedge \nz{2}{3}\wedge \nz{3}{4}\wedge \nz{3}{6}
-3\nz{0}{3}\wedge \nz{2}{3}\wedge \nz{4}{4}\wedge \nz{2}{6}
\\&
-\nz{0}{3}\wedge \nz{3}{3}\wedge \nz{0}{4}\wedge \nz{5}{6}
+2\nz{0}{3}\wedge \nz{3}{3}\wedge \nz{1}{4}\wedge \nz{4}{6}
-2\nz{0}{3}\wedge \nz{3}{3}\wedge \nz{3}{4}\wedge \nz{2}{6}
+\nz{0}{3}\wedge \nz{3}{3}\wedge \nz{4}{4}\wedge \nz{1}{6}
\\&
-3\nz{1}{3}\wedge \nz{2}{3}\wedge \nz{0}{4}\wedge \nz{5}{6}
+6\nz{1}{3}\wedge \nz{2}{3}\wedge \nz{1}{4}\wedge \nz{4}{6}
-6\nz{1}{3}\wedge \nz{2}{3}\wedge \nz{3}{4}\wedge \nz{2}{6}
+3\nz{1}{3}\wedge \nz{2}{3}\wedge \nz{4}{4}\wedge \nz{1}{6}
\\&
+3\nz{1}{3}\wedge \nz{3}{3}\wedge \nz{0}{4}\wedge \nz{4}{6}
-8\nz{1}{3}\wedge \nz{3}{3}\wedge \nz{1}{4}\wedge \nz{3}{6}
+6\nz{1}{3}\wedge \nz{3}{3}\wedge \nz{2}{4}\wedge \nz{2}{6}
-\nz{1}{3}\wedge \nz{3}{3}\wedge \nz{4}{4}\wedge \nz{0}{6}
\\&
-2\nz{2}{3}\wedge \nz{3}{3}\wedge \nz{0}{4}\wedge \nz{3}{6}
+6\nz{2}{3}\wedge \nz{3}{3}\wedge \nz{1}{4}\wedge \nz{2}{6}
-6\nz{2}{3}\wedge \nz{3}{3}\wedge \nz{2}{4}\wedge \nz{1}{6}
+2\nz{2}{3}\wedge \nz{3}{3}\wedge \nz{3}{4}\wedge \nz{0}{6}
\\
\tau_2=&3\nz{0}{3}\wedge \nz{1}{3}\wedge \nz{2}{5}\wedge \nz{5}{5}
-9\nz{0}{3}\wedge \nz{1}{3}\wedge \nz{3}{5}\wedge \nz{4}{5}
-3\nz{0}{3}\wedge \nz{2}{3}\wedge \nz{1}{5}\wedge \nz{5}{5}
+6\nz{0}{3}\wedge \nz{2}{3}\wedge \nz{2}{5}\wedge \nz{4}{5}
\\&
+\frac{2}{5}\nz{0}{3}\wedge \nz{3}{3}\wedge \nz{0}{5}\wedge \nz{5}{5}
+\nz{0}{3}\wedge \nz{3}{3}\wedge \nz{1}{5}\wedge \nz{4}{5}
-5\nz{0}{3}\wedge \nz{3}{3}\wedge \nz{2}{5}\wedge \nz{3}{5}
+\frac{9}{5}\nz{1}{3}\wedge \nz{2}{3}\wedge \nz{0}{5}\wedge \nz{5}{5}
\\&
-9\nz{1}{3}\wedge \nz{2}{3}\wedge \nz{2}{5}\wedge \nz{3}{5}
-3\nz{1}{3}\wedge \nz{3}{3}\wedge \nz{0}{5}\wedge \nz{4}{5}
+6\nz{1}{3}\wedge \nz{3}{3}\wedge \nz{1}{5}\wedge \nz{3}{5}
+3\nz{2}{3}\wedge \nz{3}{3}\wedge \nz{0}{5}\wedge \nz{3}{5}
\\&
-9\nz{2}{3}\wedge \nz{3}{3}\wedge \nz{1}{5}\wedge \nz{2}{5}
\\
\tau_3=&\nz{0}{3}\wedge \nz{1}{3}\wedge \nz{2}{5}\wedge \nz{5}{5}
-3\nz{0}{3}\wedge \nz{1}{3}\wedge \nz{3}{5}\wedge \nz{4}{5}
-\nz{0}{3}\wedge \nz{2}{3}\wedge \nz{1}{5}\wedge \nz{5}{5}
+2\nz{0}{3}\wedge \nz{2}{3}\wedge \nz{2}{5}\wedge \nz{4}{5}
\\&
+\frac{1}{5}\nz{0}{3}\wedge \nz{3}{3}\wedge \nz{0}{5}\wedge \nz{5}{5}
-\nz{0}{3}\wedge \nz{3}{3}\wedge \nz{2}{5}\wedge \nz{3}{5}
+\frac{2}{5}\nz{1}{3}\wedge \nz{2}{3}\wedge \nz{0}{5}\wedge \nz{5}{5}
+\nz{1}{3}\wedge \nz{2}{3}\wedge \nz{1}{5}\wedge \nz{4}{5}
\\&
-5\nz{1}{3}\wedge \nz{2}{3}\wedge \nz{2}{5}\wedge \nz{3}{5}
-\nz{1}{3}\wedge \nz{3}{3}\wedge \nz{0}{5}\wedge \nz{4}{5}
+2\nz{1}{3}\wedge \nz{3}{3}\wedge \nz{1}{5}\wedge \nz{3}{5}
+\nz{2}{3}\wedge \nz{3}{3}\wedge \nz{0}{5}\wedge \nz{3}{5}
\\&
-3\nz{2}{3}\wedge \nz{3}{3}\wedge \nz{1}{5}\wedge \nz{2}{5}
\end{align*}
\begin{align*} 
\tau_4=& -\frac{2}{3}\nz{0}{3}\wedge \nz{0}{4}\wedge \nz{3}{4}\wedge \nz{5}{5}
+\frac{2}{3}\nz{0}{3}\wedge \nz{0}{4}\wedge \nz{4}{4}\wedge \nz{4}{5}
+\nz{0}{3}\wedge \nz{1}{4}\wedge \nz{2}{4}\wedge \nz{5}{5}
+\frac{2}{3}\nz{0}{3}\wedge \nz{1}{4}\wedge \nz{3}{4}\wedge \nz{4}{5}
\\&
-\frac{5}{3}\nz{0}{3}\wedge \nz{1}{4}\wedge \nz{4}{4}\wedge \nz{3}{5}
-\nz{0}{3}\wedge \nz{2}{4}\wedge \nz{3}{4}\wedge \nz{3}{5}
+2\nz{0}{3}\wedge \nz{2}{4}\wedge \nz{4}{4}\wedge \nz{2}{5}
-\nz{0}{3}\wedge \nz{3}{4}\wedge \nz{4}{4}\wedge \nz{1}{5}
\\&
+\nz{1}{3}\wedge \nz{0}{4}\wedge \nz{2}{4}\wedge \nz{5}{5}
-\nz{1}{3}\wedge \nz{0}{4}\wedge \nz{4}{4}\wedge \nz{3}{5}
-7\nz{1}{3}\wedge \nz{1}{4}\wedge \nz{2}{4}\wedge \nz{4}{5}
+6\nz{1}{3}\wedge \nz{1}{4}\wedge \nz{3}{4}\wedge \nz{3}{5}
\\&
+\nz{1}{3}\wedge \nz{1}{4}\wedge \nz{4}{4}\wedge \nz{2}{5}
-\nz{1}{3}\wedge \nz{2}{4}\wedge \nz{4}{4}\wedge \nz{1}{5}
-5\nz{1}{3}\wedge \nz{2}{4}\wedge \nz{3}{4}\wedge \nz{2}{5} 
+\nz{1}{3}\wedge \nz{3}{4}\wedge \nz{4}{4}\wedge \nz{0}{5}
\\&
-\nz{2}{3}\wedge \nz{0}{4}\wedge \nz{1}{4}\wedge \nz{5}{5}
+\nz{2}{3}\wedge \nz{0}{4}\wedge \nz{2}{4}\wedge \nz{4}{5}
+\nz{2}{3}\wedge \nz{0}{4}\wedge \nz{4}{4}\wedge \nz{2}{5}
-\nz{2}{3}\wedge \nz{0}{4}\wedge \nz{3}{4}\wedge \nz{3}{5}
\\&
+5\nz{2}{3}\wedge \nz{1}{4}\wedge \nz{2}{4}\wedge \nz{3}{5}
-6\nz{2}{3}\wedge \nz{1}{4}\wedge \nz{3}{4}\wedge \nz{2}{5}
+7\nz{2}{3}\wedge \nz{2}{4}\wedge \nz{3}{4}\wedge \nz{1}{5}
-\nz{2}{3}\wedge \nz{2}{4}\wedge \nz{4}{4}\wedge \nz{0}{5}
\\&
+\nz{3}{3}\wedge \nz{0}{4}\wedge \nz{1}{4}\wedge \nz{4}{5}
-2\nz{3}{3}\wedge \nz{0}{4}\wedge \nz{2}{4}\wedge \nz{3}{5}
+\frac{5}{3}\nz{3}{3}\wedge \nz{0}{4}\wedge \nz{3}{4}\wedge \nz{2}{5}
-\frac{2}{3}\nz{3}{3}\wedge \nz{0}{4}\wedge \nz{4}{4}\wedge \nz{1}{5}
\\&
+\nz{3}{3}\wedge \nz{1}{4}\wedge \nz{2}{4}\wedge \nz{2}{5}
-\frac{2}{3}\nz{3}{3}\wedge \nz{1}{4}\wedge \nz{3}{4}\wedge \nz{1}{5}
+\frac{2}{3}\nz{3}{3}\wedge \nz{1}{4}\wedge \nz{4}{4}\wedge \nz{0}{5}
-\nz{3}{3}\wedge \nz{2}{4}\wedge \nz{3}{4}\wedge \nz{0}{5}
\\
\tau_5=&\nz{0}{3}\wedge \nz{0}{4}\wedge \nz{3}{4}\wedge \nz{5}{5}
-\nz{0}{3}\wedge \nz{0}{4}\wedge \nz{4}{4}\wedge \nz{4}{5}
-4\nz{0}{3}\wedge \nz{1}{4}\wedge \nz{3}{4}\wedge \nz{4}{5}
+4\nz{0}{3}\wedge \nz{1}{4}\wedge \nz{4}{4}\wedge \nz{3}{5}
\\&
+6\nz{0}{3}\wedge \nz{2}{4}\wedge \nz{3}{4}\wedge \nz{3}{5}
-6\nz{0}{3}\wedge \nz{2}{4}\wedge \nz{4}{4}\wedge \nz{2}{5}
+3\nz{0}{3}\wedge \nz{3}{4}\wedge \nz{4}{4}\wedge \nz{1}{5}
-3\nz{1}{3}\wedge \nz{0}{4}\wedge \nz{2}{4}\wedge \nz{5}{5}
\\&
+12\nz{1}{3}\wedge \nz{1}{4}\wedge \nz{2}{4}\wedge \nz{4}{5}
-12\nz{1}{3}\wedge \nz{1}{4}\wedge \nz{3}{4}\wedge \nz{3}{5}
+3\nz{1}{3}\wedge \nz{2}{4}\wedge \nz{4}{4}\wedge \nz{1}{5}
+6\nz{1}{3}\wedge \nz{2}{4}\wedge \nz{3}{4}\wedge \nz{2}{5}
\\&
-3\nz{1}{3}\wedge \nz{3}{4}\wedge \nz{4}{4}\wedge \nz{0}{5}
+3\nz{2}{3}\wedge \nz{0}{4}\wedge \nz{1}{4}\wedge \nz{5}{5}
-3\nz{2}{3}\wedge \nz{0}{4}\wedge \nz{2}{4}\wedge \nz{4}{5}
-6\nz{2}{3}\wedge \nz{1}{4}\wedge \nz{2}{4}\wedge \nz{3}{5}
\\&
+12\nz{2}{3}\wedge \nz{1}{4}\wedge \nz{3}{4}\wedge \nz{2}{5}
-12\nz{2}{3}\wedge \nz{2}{4}\wedge \nz{3}{4}\wedge \nz{1}{5}
+3\nz{2}{3}\wedge \nz{2}{4}\wedge \nz{4}{4}\wedge \nz{0}{5}
-3\nz{3}{3}\wedge \nz{0}{4}\wedge \nz{1}{4}\wedge \nz{4}{5}
\\&
+6\nz{3}{3}\wedge \nz{0}{4}\wedge \nz{2}{4}\wedge \nz{3}{5}
-4\nz{3}{3}\wedge \nz{0}{4}\wedge \nz{3}{4}\wedge \nz{2}{5}
+\nz{3}{3}\wedge \nz{0}{4}\wedge \nz{4}{4}\wedge \nz{1}{5}
-6\nz{3}{3}\wedge \nz{1}{4}\wedge \nz{2}{4}\wedge \nz{2}{5}
\\&
+4\nz{3}{3}\wedge \nz{1}{4}\wedge \nz{3}{4}\wedge \nz{1}{5}
-\nz{3}{3}\wedge \nz{1}{4}\wedge \nz{4}{4}\wedge \nz{0}{5}
+3\nz{1}{3}\wedge \nz{0}{4}\wedge \nz{3}{4}\wedge \nz{4}{5}
-3\nz{2}{3}\wedge \nz{1}{4}\wedge \nz{4}{4}\wedge \nz{1}{5}
\end{align*}

%% file: in_dp5_w8.tex
\begin{align*}
\rho =&3\nz{0}{3}\wedge \nz{1}{3}\wedge \nz{2}{3}\wedge \nz{1}{4}\wedge \nz{5}{5}
-9\nz{0}{3}\wedge \nz{1}{3}\wedge \nz{2}{3}\wedge \nz{2}{4}\wedge \nz{4}{5}
+9\nz{0}{3}\wedge \nz{1}{3}\wedge \nz{2}{3}\wedge \nz{3}{4}\wedge \nz{3}{5}
\\&
-3\nz{0}{3}\wedge \nz{1}{3}\wedge \nz{2}{3}\wedge \nz{4}{4}\wedge \nz{2}{5}
-\nz{0}{3}\wedge \nz{1}{3}\wedge \nz{3}{3}\wedge \nz{0}{4}\wedge \nz{5}{5}
+\nz{0}{3}\wedge \nz{1}{3}\wedge \nz{3}{3}\wedge \nz{1}{4}\wedge \nz{4}{5}
\\&
+3\nz{0}{3}\wedge \nz{1}{3}\wedge \nz{3}{3}\wedge \nz{2}{4}\wedge \nz{3}{5}
-5\nz{0}{3}\wedge \nz{1}{3}\wedge \nz{3}{3}\wedge \nz{3}{4}\wedge \nz{2}{5}
+2\nz{0}{3}\wedge \nz{1}{3}\wedge \nz{3}{3}\wedge \nz{4}{4}\wedge \nz{1}{5}
\\&
+2\nz{0}{3}\wedge \nz{2}{3}\wedge \nz{3}{3}\wedge \nz{0}{4}\wedge \nz{4}{5}
-5\nz{0}{3}\wedge \nz{2}{3}\wedge \nz{3}{3}\wedge \nz{1}{4}\wedge \nz{3}{5}
+3\nz{0}{3}\wedge \nz{2}{3}\wedge \nz{3}{3}\wedge \nz{2}{4}\wedge \nz{2}{5}
\\&
+\nz{0}{3}\wedge \nz{2}{3}\wedge \nz{3}{3}\wedge \nz{3}{4}\wedge \nz{1}{5}
-\nz{0}{3}\wedge \nz{2}{3}\wedge \nz{3}{3}\wedge \nz{4}{4}\wedge \nz{0}{5}
-3\nz{1}{3}\wedge \nz{2}{3}\wedge \nz{3}{3}\wedge \nz{0}{4}\wedge \nz{3}{5}
\\&
+9\nz{1}{3}\wedge \nz{2}{3}\wedge \nz{3}{3}\wedge \nz{1}{4}\wedge \nz{2}{5}
-9\nz{1}{3}\wedge \nz{2}{3}\wedge \nz{3}{3}\wedge \nz{2}{4}\wedge \nz{1}{5}
+3\nz{1}{3}\wedge \nz{2}{3}\wedge \nz{3}{3}\wedge \nz{3}{4}\wedge \nz{0}{5}
\end{align*}